\documentclass[smallcondensed]{svjour3}
\usepackage{amssymb}
\usepackage[dvips]{graphics}
\input{amssym.def}
\input{amssym.tex}
\title{\huge\bf Curvilinear Mode-I/Mode-II interface fracture with a curvature-dependent surface tension on the boundary}
\author{Anna Y. Zemlyanova}
\institute{Department of Mathematics, Kansas State University, 138 Cardwell Hall, Manhattan KS 66506\\
Tel.: +1-785-532-6750, Fax: +1-785-532-0546\\
\email{azem@math.ksu.edu}}

\titlerunning{Curvilinear Mode-I/Mode-II interface fracture with a surface tension}
\authorrunning{Anna Y. Zemlyanova}
\journalname{}

\date{}

\newcommand{\I}{\mathop{\rm Im}\nolimits}
\newcommand{\R}{\mathop{\rm Re}\nolimits}

\newcommand{\bfa}[1]{\mathbf{#1}}     %
\newcommand{\grad}{\mbox{grad}}       %
\newcommand{\divv}{\mbox{div}}        %
\newcommand{\tgamma}{{\tilde{\gamma}}}%
\newcommand{\jump}[1]{[\![ #1 ]\!]}

\begin{document}
\maketitle

\noindent

\begin{abstract}
A new model of fracture mechanics considered previously by Sendova and Walton \cite{SendovaWalton2010}, Zemlyanova \cite{Zemlyanova2013}, and Zemlyanova and Walton \cite{Zemlyanova2012} is further developed on the example of a mixed mode curvilinear interface fracture located on the boundary of a partially debonded thin elastic inclusion embedded in an infinite thin elastic matrix. The effect of the nano-structure of the material is incorporated into the model in the form of a curvature-depended surface tension acting on the boundary of the fracture. It is shown that the introduction of the surface tension allows to eliminate the classical oscillating and power singularities of the order $1/2$ present in the linear elastic fracture mechanics. The mathematical methods used to solve the problem are based on the Muskhelishvili's complex potentials and the Savruk's integral representations. The mechanical problem is reduced to the system of singular integro-differential equations which is further reduced to a system of weakly-singular integral equations. The numerical computations and comparison with known results are presented.

\end{abstract}



\setcounter{equation}{0}

\section{Introduction}

The study of brittle fracture in solids has been a subject of many investigations. Historically, the fracture mechanics problems were considered within the framework of linear elastic fracture mechanics (LEFM). The LEFM theory allows to reduce many mechanical problems to relatively simple mathematical equations which can be studied efficiently both theoretically and numerically while describing the behavior of many materials with sufficient accuracy. At the same time, LEFM contains important internal inconsistencies. In particular, while all of the equations of linear elasticity are obtained under assumption that the stresses and the strains are bounded everywhere in an object, LEFM predicts a power singularity of the order $1/2$ in the stresses and the strains at the crack tips. This contradiction becomes even more pronounced for an interface fracture on the boundary of two materials with different mechanical properties. In this case, in addition to the power singularity of the order $1/2$, the stresses and the strains possess a power singularity of a pure imaginary order (oscillating singularity) which predicts wrinkling and interpenetration of the crack surfaces near the crack tips.

Multiple attempts have been made to remove these inconsistencies of LEFM. It has been observed that due to the stress concentration near the crack tips, the behavior of the material can no longer be considered linearly elastic. Consequently, cohesive and processing zone models have been introduced to account for non-linear elastic and plastic deformations as well as damage accumulation near the crack tips. The main difficulty in the application of these models is in specification of the constitutive properties in the cohesive or processing zones which are very difficult to obtain from experiment. Extensive work in this area has been done by many authors including classical papers by Barenblatt \cite{Barenblatt1962} and Dugdale \cite{Dugdale1960}.

Since the fracture initiation and propagation is a nano-scale process, it has been argued that continuum models cannot be used to accurately describe it. Various atomistic and lattice approaches have been developed to overcome this difficulty \cite{Abraham2001}, \cite{FinebergGross1991}, \cite{HollandMarder1998}, \cite{SlepyanEtAl1999}. This approach requires an accurate description of long-range and short-range intermolecular forces  and also presents some computational challenges.

Blended atom-to-continuum models have recently generated a considerable attention.  Among those it is necessary to point out the quasi-continuum method introduced by Tadmor et al \cite{TadmoretAl1996} and the method based on the introduction of bridging domains between continuum and atomistic regions proposed by Xiao and Belytschko \cite{XiaoBelytschko2004}. Atom-to-continuum models involve adjustable parameters,  such as choosing the size and location of domains over which  potentials acting at different length scales are blended. One of the main difficulties here is in constructing a coupled energy whose minimizers are free from uncontrollable errors (``ghost forces") on the atomistic/continuum interface. One of the ways to overcome this difficulty has been proposed in \cite{Shapeev2011}.

Oscillating power singularity of a pure imaginary order at the crack tips of an interface fracture has been first described by Williams in \cite{Williams1959}, and later confirmed by other authors. Different approaches have been proposed to eliminate this particular type of singularity such as the contact zone model suggested by Comninou \cite{Comninou1990} and the intermediate layer model suggested by Atkinson \cite{Atkinson1977}.

Recently, several continuum models of fracture mechanics with surface excess properties on the boundary of the fracture have been studied as well. The physical motivation behind these models stems from the fact that the material particles on the boundary of the solid experience different force system compared to the particles in the bulk. The first comprehensive model of surface stressed solids has been developed by Gurtin and Murdoch \cite{GurtinMurdoch1975}, \cite{GurtinMurdoch1978}. Their approach is based on the notion of a thin two-dimensional membrane that is effectively bonded to the surface of a three-dimensional bulk substrate. The membrane is modeled as an elastic surface without accounting explicitly for its thickness in three dimensions. This idea generalizes the classical notions of a fracture energy and surface tension and furnishes a basis for study of mechanics of coated surfaces. 

The Gurtin-Murdoch theory is very popular and has been widely applied to the study of inhomogeneities and nano-structures by many authors \cite{Duan2005a}, \cite{Mogilevskaya2010}, \cite{Sharma2004}, \cite{Tian2007}. Recently Kim, Schiavone, and Ru  applied the Gurtin-Murdoch theory to the modeling of fracture \cite{Schiavone2010a}-\cite{Schiavone2012}. Their results for a straight interface and a non-interface mode-III fracture show that taking into account the surface elasticity eliminates the square root singularities of the stresses and the strains. At the same time it has been shown \cite{Schiavone2012} that the presence of the surface elasticity by itself is not sufficient to guarantee that the stresses and the strains are bounded at the crack tips. Weaker logarithmic singularities may still be present. Moreover, a mode-I fracture may still continue to exhibit strong square-root singularities even if the surface elasticity is taken into account \cite{Schiavone2012}.

Another theory of solids with surface excess properties has been developed by Slattery et al \cite{SlatteryEtAl2004}, \cite{Slatteryetal2007}. It is assumed that the material surfaces are endowed with a curvature-dependent surface tension. This theory has been first applied to the study of fracture in \cite{OhWaltonSlattery2006} and developed in more detail on the example of a straight mode-I fracture in \cite{SendovaWalton2010}. The results of Sendova and Walton \cite{SendovaWalton2010} show that for a mode-I straight non-interface fracture all of the stresses and strains are bounded if the surface tension on the boundary of the fracture depends linearly on the mean curvature of the deformed fracture. However, it has been shown later  on the example of a mode-I/mode-II non-interface curvilinear plane fracture \cite{Zemlyanova2012} and on the example of a straight mode-I/mode-II interface fracture  \cite{Zemlyanova2013} that this conclusion does not hold for non-symmetric configurations. The main conclusion of the papers \cite{SendovaWalton2010}, \cite{Zemlyanova2013}, \cite{Zemlyanova2012} is  that incorporation of the curvature-dependent surface tension allows to eliminate the classical singularity of the order $1/2$ and also oscillating singularity in the case of a straight interface crack. However, some components of the stresses and strains may still retain a weaker logarithmic singularity. It has been shown recently \cite{Walton2013} that the modification of the curvature-dependent surface tension model, which includes a surface stretch in addition to the mean curvature, removes the logarithmic singularities in the case of a straight mode-I/mode-II fracture.

The determination of the surface energy parameters presents an interesting practical problem. One of a few studies dedicated to this problem is the paper by Mohammadi and Sharma \cite{Mohammadi2012}. The surface energy parameters have been obtained within the framework of the Steigmann-Ogden theory by comparing theoretical results with atomistic computations for nanowires with flat and corrugated surfaces. It has been shown in \cite{Mohammadi2012} that surface roughness significantly influences the values of the surface energy parameters. Additional discussion of these issues is available in \cite{Ergincanetal} and \cite{WeisDuan2008}.

The focus of the current study is to continue the investigation started in \cite{Zemlyanova2013}, \cite{Zemlyanova2012} and to  generalize the curvature-depended surface tension fracture model proposed in \cite{OhWaltonSlattery2006}, \cite{SendovaWalton2010} to the case of an interface curvilinear fracture of an arbitrary shape on the boundary of a partially debonded thin elastic inclusion embedded into a thin elastic plate (matrix). It is assumed that in-plane stresses act on the boundary of the fracture and at infinity of the plate. The nano-structure of the material of the inclusion and the matrix is modeled with the help of a curvature-depended surface tension which acts on the interface between two different materials or between a material and a void. It is assumed that the fracture can be of arbitrary smooth shape. This is a significant advantage of the present study since most of the fracture literature deals with straight or circular-arc-shaped cracks. The mechanical problem is reduced to a system of singular integro-differential equations using methods of complex analysis. This system is further reduced to a system of weakly-singular integral equations. The numerical solution to the system of singular integro-differential equations is obtained using approximations of the unknown functions by Taylor polynomials with unknown coefficients. This allows to reduce the approximate solution of the system of integro-differential equations to the system of linear algebraic equations. Convergence of the proposed numerical method is studied on examples. Comparison of the current results with known results obtained in \cite{CrouchMogilevskaya2006}, \cite{PrasadSimha2003}, \cite{SilZem2001}, \cite{Toya1974}, \cite{Zemlyanova2012} is given as well. The main conclusion of the paper is that the introduction of the curvature-dependent surface tension allows to eliminate integrable power singularities of the order $1/2$ and oscillating singularities of the stresses and the strains at the crack-tips, and hence presents a significant improvement of the classical LEFM theory.

\setcounter{equation}{0}

\section{Model with a curvature-dependent surface tension}

\begin{figure}[ht]
	\centering
		\scalebox{0.5}{\includegraphics{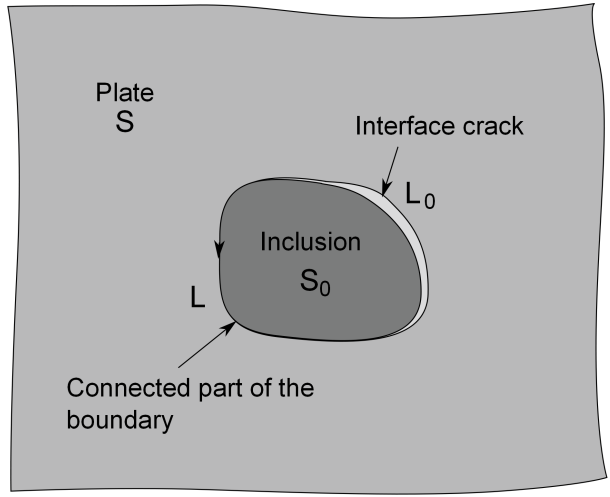}}
			\caption{An infinite plate $S$ with an inclusion $S_0$ and an interface crack $L_0$}
	\label{fig0}
\end{figure}

Consider an infinite thin elastic plate (matrix) $S$ which has a hole with a smooth boundary $L_0\cup L$. A thin elastic inclusion $S_0$ of the same size as the hole is inserted into the hole and partially attached to the matrix along the line $L$ (fig.~\ref{fig0}). The plate and the inclusion are homogeneous and isotropic. Their shear moduli and Poisson ratios are given by $\mu$, $\nu$ and $\mu_0$, $\nu_0$ correspondingly. The principal in-plane stresses $\sigma_1^{\infty}$ and $\sigma_2^{\infty}$ are applied at infinity of the plate and act in the directions constituting the angles $\alpha$ and $\alpha+\pi/2$ with the positive direction of the real axis. 

Assume that the inclusion $S_0$ and the matrix $S$ are perfectly attached along the junction line $L$:
\begin{equation}
(u_1+iu_2)_0^{+}(t)=(u_1+iu_2)^{-}(t), \, \, t\in L,
\label{2_1}
\end{equation}    
where $(u_1+iu_2)(t)$ is the vector of the displacements at a point $t$ of the plate $S$ or the inclusion $S_0$. Parameters without a subscript are related to the plate $S$; parameters with a subscript ``$0$" are related to the inclusion $S_0$. Here and henceforth, the superscripts ``$+$" and ``$-$" denote the limit values of the stresses, the displacements and other parameters from the left-hand and the right-hand sides correspondingly of the lines $L_0$ or $L$. The default direction of the curve $L_0\cup L$ is chosen to be counterclockwise (fig.~\ref{fig0}).

Assume that the surface tension $\tilde{\gamma}$ acts on the dividing lines $L$ and $L_0$ which separate materials with different properties (two different elastic materials on the line $L$ and an elastic material and a void/gas on the left- or the right-hand side of the line $L_0$). The surface tension $\tilde{\gamma}$ allows to take into the account the effects of long-range intermolecular forces on the dividing lines $L_0$, $L$. Then the differential and the jump momentum conditions in the deformed configuration in the absence of inertial and gravitational effects become \cite{Slatteryetal2007}:
\begin{equation}
\mbox{div}(\bfa{T})=0,
\label{2_2}
\end{equation}
\begin{equation}
\grad_{(\zeta)}\tgamma+2\tgamma H\bfa{n} +\jump{\bfa{T}}\bfa{n} =\bfa{0},
\label{2_3}
\end{equation}
where $\bfa{T}$ is the Cauchy stress tensor, $\bfa{n}$ is the unit normal to the fracture surface $\zeta$ pointing into the bulk of the material, $H=-\frac12 \mbox{div}_{(\zeta)}\bfa{n}$ is the mean curvature, $\grad_{(\zeta)}$, $\mbox{div}_{(\zeta)}$ denote the surface gradient and the surface divergence correspondingly, and the double brackets $\jump{\ldots}$ denote the jump of the quantity enclosed across the boundary of the line $L_0\cup L$. The equation (\ref{2_2}) is valid in the bulk of the material of the inclusion $S_0$ or the matrix $S$, and the equation (\ref{2_3}) is valid on the boundary lines $L_0$ and $L$.

Assume that the surface tension $\tgamma$ depends linearly on the difference between the curvature of the deformed line $L_0\cup L$ and the curvature of the line $L_0\cup L$ in the unloaded configuration \cite{Zemlyanova2012}:

\begin{equation}
\tilde{\gamma}=\gamma^{*}(\divv_{(\zeta)}{\bfa{n}}-\divv_{(\zeta_{in})}{\bfa{n}_{in}}),\,\,\,\gamma^{*}=\mbox{const},
\label{2_3a}
\end{equation}
where the subindex ``$in$" denotes the parameters in the initial undeformed configuration of the line $L_0\cup L$, and the superindex ``$*$" should be replaced by ``$+$" to denote the surface tension parameter on the crack $L_0$ from the side of the inclusion $S_0$, by ``$-$" to denote the surface tension parameter on the crack $L_0$ from the side of the plate $S$, and by ``$i$" to denote the surface tension parameter between the inclusion $S_0$ and the matrix $S$ on the junction line $L$.

Let $t=t(s)$ be the parametric equation of the line $L_0\cup L$ (fig. \ref{fig0}), where the parameter $s$ is an arc length, and the values $s\in [0,l_0]$ correspond to the curve $L_0$, while the values $s\in [l_0, l]$ correspond to the curve $L$. Since the contour $L_0\cup L$ is closed, it is necessary that $t(0)=t(l)$. Assume that the function $t(s)$ has continuous derivatives up to the fourth order. This assumption is made for simplicity and can be somewhat relaxed. It is possible to rewrite the condition (\ref{2_3}) in terms of the normal and shear stresses in the following form:
\begin{equation}
\sigma_n+i\tau_n=i(\grad_{(\zeta)}\tgamma-\tgamma\bfa{n}\divv_{(\zeta)}\bfa{n})\overline{t'(s)},
\label{2_4}
\end{equation}
where $\sigma_n$ and $\tau_n$ are the tensile and the shear components of the stress vector acting on the tangent line to the curves $L_0$ and $L$.

Linearizing the condition (\ref{2_4}) under the assumption that the derivatives of the displacements are small, similarly to \cite{Zemlyanova2012}, obtain the following boundary conditions on the lines $L_0$ and $L$:
$$
(\sigma_n+i\tau_n)_0^{+}(s)=\frac{\gamma^+}{2}m_1(s)\frac{d}{ds}(u_1+iu_2)_0^{+}(s)+\frac{\gamma^+}{2}m_2(s)\frac{d}{ds}{(u_1-iu_2)_0^{+}(s)}+
$$
$$
\frac{\gamma^+}{2}m_3(s)\frac{d^2}{ds^2}(u_1+iu_2)_0^{+}(s)+\frac{\gamma^+}{2}m_4(s)\frac{d^2}{ds^2}{(u_1-iu_2)_0^{+}(s)}+
$$
\begin{equation}
\frac{\gamma^+}{2}\overline{t'(s)}\frac{d^3}{ds^3}(u_1+iu_2)_0^{+}(s)-\frac{\gamma^+}{2}t'(s)\frac{d^3}{ds^3}{(u_1-iu_2)_0^{+}(s)}+f_1(s),\,\,\,\,s\in [0,l_0],
\label{2_5}
\end{equation}
$$
(\sigma_n+i\tau_n)^{-}(s)=\frac{\gamma^-}{2}m_1(s)\frac{d}{ds}(u_1+iu_2)^{-}(s)+\frac{\gamma^-}{2}m_2(s)\frac{d}{ds}{(u_1-iu_2)}^{-}(s)+
$$
$$
\frac{\gamma^-}{2}m_3(s)\frac{d^2}{ds^2}(u_1+iu_2)^{-}(s)+\frac{\gamma^-}{2}m_4(s)\frac{d^2}{ds^2}{(u_1-iu_2)}^{-}(s)+
$$
\begin{equation}
\frac{\gamma^-}{2}\overline{t'(s)}\frac{d^3}{ds^3}(u_1+iu_2)^{-}(s)-\frac{\gamma^-}{2}t'(s)\frac{d^3}{ds^3}{(u_1-iu_2)}^{-}(s)+f_2(s),\,\,\,\,s\in [0,l_0],
\label{2_6}
\end{equation}
$$
(\sigma_n+i\tau_n)_0^{+}(s)-(\sigma_n+i\tau_n)^{-}(s)=\frac{\gamma^i}{2}m_1(s)\frac{d}{ds}(u_1+iu_2)_0^{+}(s)+
$$
$$
\frac{\gamma^i}{2}m_2(s)\frac{d}{ds}{(u_1-iu_2)_0^{+}(s)}+\frac{\gamma^i}{2}m_3(s)\frac{d^2}{ds^2}(u_1+iu_2)_0^{+}(s)+
$$
\begin{equation}
\frac{\gamma^i}{2}m_4(s)\frac{d^2}{ds^2}{(u_1-iu_2)_0^{+}(s)}+\frac{\gamma^i}{2}\overline{t'(s)}\frac{d^3}{ds^3}(u_1+iu_2)_0^{+}(s)-
\label{2_7}
\end{equation}
$$
\frac{\gamma^i}{2}t'(s)\frac{d^3}{ds^3}{(u_1-iu_2)_0^{+}(s)},\,\,\,\,s\in [l_0,l],
$$
where
$$
m_1(s)=-\overline{t'''(s)}-2i\overline{t''(s)}\varrho(s)-3i\overline{t'(s)}\varrho'(s)-3\overline{t'(s)}\varrho^2(s),
$$
$$
m_2(s)={t'''(s)}-4i{t''(s)}\varrho(s)-3i{t'(s)}\varrho'(s)-3{t'(s)}\varrho^2(s),
$$
$$
m_3(s)=-4i\overline{t'(s)}\varrho(s),\,\,\,\,m_4(s)=-2i{t'(s)}\varrho(s),
$$
$f_1(s)$ and $f_2(s)$ are given external in-plane stresses applied to the banks of the crack, the superindices ``$+$", ``$-$", and ``$i$" have the same meaning as in (\ref{2_3a}), $\varrho(s)$ denotes the mean curvature of the  line $L_0$ or $L$ in the undeformed initial configuration at the point described by the arc length $s$. The equations (\ref{2_5}) and (\ref{2_6}) correspond to the boundary conditions on the crack $L_0$ from the side of the inclusion $S_0$ and the matrix $S$ correspondingly, the equation (\ref{2_7}) describes the boundary condition on the junction line $L$.

\setcounter{equation}{0}
\section{Complex potentials and their integral representations}

The presented method of the solution is based on the use of the Muskhelishvili's complex potentials \cite{Mus1963} and the Savruk's integral representations \cite{Savruk1981}. It is well-known that the stresses $\sigma_n$ and $\tau_n$ and the derivatives of the displacements $d(u_1+iu_2)/dt$ in the inclusion $S_0$ and the matrix $S$ can be expressed through two complex functions $\Phi(z)$, $\Psi(z)$ (complex potentials) analytic in $S_0$ or $S$ correspondingly using the following formulas \cite{Mus1963}:\\
\begin{equation}
\label{3_1}
(\sigma_n+i\tau_n)(t)=\Phi(t)+\overline{\Phi(t)}+\frac{\overline{dt}}{dt}(t\overline{\Phi'(t)}+\overline{\Psi(t)}),
\end{equation}
\begin{equation}
\label{3_2}
2\mu\frac{d}{dt}(u_1+iu_2)(t)=\kappa\Phi(t)-\overline{\Phi(t)}-\frac{\overline{dt}}{dt}(t\overline{\Phi'(t)}+\overline{\Psi(t)}), \,\, t\in L_0\cup L.
\end{equation}
Here $\kappa=(3-\nu)/(1+\nu)$ for the case of plane stress and $\kappa=3-4\nu$ for the plane strain, and all of the parameters and functions corresponding to the inclusion $S_0$ should be taken with ``$0$" subscript, for example, $(\sigma_n+i\tau_n)_0$, $(u_1+iu_2)_0$, $\Phi_0(t)$, $\Psi_0(t)$, $\mu_0$, $\kappa_0$ and so on.  

The analytic functions $\Phi_0(z)$, $\Psi_0(z)$ and $\Phi(z)$, $\Psi(z)$ can be expressed through the jumps of the stresses and the derivatives of the displacements on the contour $L_0\cup L$ by using the following integral representations which are due to Savruk \cite{Savruk1981}: 
\begin{equation}
\label{3_3}
\Phi_0(z)=\frac{1}{2\pi}\int_{L_0\cup L}\frac{g'_0(t)dt}{t-z}+\frac{(\kappa_0+1)^{-1}}{\pi i}\int_{L_0\cup L}\frac{q_0(t)dt}{t-z},
\end{equation}
$$
\Psi_0(z)=\frac{1}{2\pi}\int_{L_0\cup L}\left( \frac{\overline{g'_0(t)dt}}{t-z}-\frac{\bar{t}g'_0(t)dt}{(t-z)^2}\right)+
$$
$$
\frac{(\kappa_0+1)^{-1}}{\pi i}\int_{L_0\cup L}\left(\frac{\kappa_0\overline{q_0(t)dt}}{t-z}-\frac{\bar{t}q_0(t)dt}{(t-z)^2}\right), \,\, z\in S_0,
$$
\begin{equation}
\label{3_4}
\Phi(z)=\Gamma+\frac{1}{2\pi}\int_{L_0\cup L}\frac{g'(t)dt}{t-z}+\frac{(\kappa+1)^{-1}}{\pi i}\int_{L_0\cup L}\frac{q(t)dt}{t-z},
\end{equation}
$$
\Psi(z)=\Gamma'+\frac{1}{2\pi}\int_{L_0\cup L}\left( \frac{\overline{g'(t)dt}}{t-z}-\frac{\bar{t}g'(t)dt}{(t-z)^2}\right)+
$$
$$
\frac{(\kappa+1)^{-1}}{\pi i}\int_{L_0\cup L}\left(\frac{\kappa\overline{q(t)dt}}{t-z}-\frac{\bar{t}q(t)dt}{(t-z)^2}\right), \,\, z\in S,
$$
$$
\Gamma=(\sigma_1^{\infty}+\sigma_2^{\infty})/4, \,\, \Gamma'=(\sigma_2^{\infty}-\sigma_1^{\infty})e^{-2i\alpha}/2,
$$
where the jumps of the stresses $q_0(t)$, $q(t)$ and the derivatives of the displacements $g'_0(t)$, $g'(t)$ on the contour $L_0\cup L$ in the inclusion $S_0$ and in the plate $S$ can be found from the formulas:
\begin{equation}
2q_0(t)=(\sigma_n+i\tau_n)_0^+(t)-(\sigma_n+i\tau_n)_0^-(t),\,\,\,t\in L_0\cup L,
\label{3_5}
\end{equation}
\begin{equation}
\frac{i(\kappa_0+1)}{2\mu_0}g'_0(t)=\frac{d}{dt}(u_1+iu_2)_0^+(t)-\frac{d}{dt}(u_1+iu_2)_0^-(t),\,\,\, t\in L_0\cup L,
\label{3_6}
\end{equation}
\begin{equation}
2q(t)=(\sigma_n+i\tau_n)^+(t)-(\sigma_n+i\tau_n)^-(t),\,\,\,t\in L_0\cup L,
\label{3_7}
\end{equation}
\begin{equation}
\frac{i(\kappa+1)}{2\mu}g'(t)=\frac{d}{dt}(u_1+iu_2)^+(t)-\frac{d}{dt}(u_1+iu_2)^-(t),\,\,\, t\in L_0\cup L.
\label{3_8}
\end{equation}

Hence, the stressed state of the inclusion $S_0$ and the matrix $S$ is described by the complex potentials (\ref{3_3}), (\ref{3_4}) which contain four unknown functions $g'_0(t)$, $q_0(t)$ and $g'(t)$, $q(t)$ defined on the contour $L_0\cup L$. We will look for these functions in the class of functions satisfying the H\"older condition on the curves $L_0$, $L$. This choice guarantees the existence of all principal and limit values of the integrals of the Cauchy type in the formulas (\ref{3_3}), (\ref{3_4}), except for, maybe, at the end-points of the curves $L_0$, $L$. 

The matrix $S$ occupies the exterior of the line $L_0\cup L$, and the inclusion $S_0$ occupies the interior of the line $L_0\cup L$. Hence, the values $(\sigma_n+i\tau_n)_0^-(t)$, $(\sigma_n+i\tau_n)^+(t)$, $\frac{d}{dt}(u_1+iu_2)_0^-(t)$, $\frac{d}{dt}(u_1+iu_2)^+(t)$ in the formulas (\ref{3_5})-(\ref{3_8}) do not have any physical meaning and can be defined formally in multiple ways. The choice of the definition depends on the ease of treatment of the resulting system of singular integro-differential equations. In this case, it is convenient to formally extend the inclusion $S_0$ and the matrix $S$ to a full complex plane by assuming that the stresses and the derivatives of the displacements are equal to zero outside of the line $L_0\cup L$ for the inclusion $S_0$ or inside of the line $L_0\cup L$ for the matrix $S$ correspondingly:
\begin{equation}
(\sigma_n+i\tau_n)_0(t)=0,\,\,\,\frac{d}{dt}(u_1+iu_2)_0(t)=0,\,\,\,t\in {\mathbb{C}}\setminus S_0,
\label{3_9}
\end{equation}
$$
(\sigma_n+i\tau_n)(t)=0,\,\,\,\frac{d}{dt}(u_1+iu_2)(t)=0,\,\,\,t\in {\mathbb{C}}\setminus S.
$$
In particular, this leads to the boundary conditions: 
\begin{equation}
\frac{d}{dt}(u_1+iu_2)^-_0(t)=0,\,\,\,\frac{d}{dt}(u_1+iu_2)^+(t)=0,\,\,\,t\in L_0\cup L,
\label{3_10}
\end{equation}
\begin{equation}
(\sigma_n+i\tau_n)^-_0(t)=0,\,\,\,(\sigma_n+i\tau_n)^+(t)=0,\,\,\,t\in L_0\cup L.
\label{3_11}
\end{equation}
Observe, that if we assume the boundary conditions (\ref{3_10}), then from the uniqueness of the solution of the second fundamental problem of elasticity \cite{Mus1963}, it follows that the conditions (\ref{3_11}) are true as well and that the equations (\ref{3_9}) are satisfied.

Observe also, that substituting the equations (\ref{3_10}), (\ref{3_11}) into the relations (\ref{3_5})-(\ref{3_8}) leads to the convenient boundary conditions:
\begin{equation}
\frac{d}{dt}(u_1+iu_2)^+_0(t)=\frac{i(\kappa_0+1)}{2\mu_0}g'_0(t),\,\,\,\frac{d}{dt}(u_1+iu_2)^-(t)=-\frac{i(\kappa+1)}{2\mu}g'(t),
\label{3_10a}
\end{equation}
$$
t\in L_0\cup L,
$$
\begin{equation}
(\sigma_n+i\tau_n)^+_0(t)=2q_0(t),\,\,\,(\sigma_n+i\tau_n)^-(t)=-2q(t),\,\,\,t\in L_0\cup L.
\label{3_11a}
\end{equation}

If the equations (\ref{3_10a})-(\ref{3_11a}) are further substituted into the equations (\ref{2_5})-(\ref{2_7}) and the equation (\ref{2_1}) differentiated by $t$, the conditions on the boundaries of the inclusion $S_0$ and the matrix $S$ become:
\begin{equation}
\R q_0(s)=\frac{\gamma^+}{4\mu_0}(\kappa_0+1)\varrho(s)\left[\varrho(s)\I g'_0(s)+\R g''_0(s) \right]+
\label{3_12a}
\end{equation}
$$
\frac12\R f_1(s),\,\,s\in [0,l_0],
$$
\begin{equation}
\I q_0(s)=\frac{\gamma^+}{4\mu_0}(\kappa_0+1)\frac{d}{ds}\left[\varrho(s)\I g'_0(s)+\R g''_0(s) \right]+\frac12\I f_1(s),\,\,s\in [0,l_0],
\label{3_12}
\end{equation}
\begin{equation}
\R q(s)=\frac{\gamma^-}{4\mu}(\kappa+1)\varrho(s)\left[\varrho(s)\I g'(s)+\R g''(s) \right]-\frac12\R f_2(s),\,\,s\in [0,l_0],
\label{3_13a}
\end{equation}
\begin{equation}
\I q(s)=\frac{\gamma^-}{4\mu}(\kappa+1)\frac{d}{ds}\left[\varrho(s)\I g'(s)+\R g''(s) \right]-\frac12\I f_2(s),\,\,s\in [0,l_0],
\label{3_13}
\end{equation}
\begin{equation}
\R q_0(s)+\R q(s)=\frac{\gamma^i}{4\mu_0}(\kappa_0+1)\varrho(s)\left[\varrho(s)\I g'_0(s)+\R g''_0(s) \right],\,\,s\in [l_0,l],
\label{3_14a}
\end{equation}
\begin{equation}
\I q_0(s)+\I q(s)=\frac{\gamma^i}{4\mu_0}(\kappa_0+1)\frac{d}{ds}\left[\varrho(s)\I g'_0(s)+\R g''_0(s) \right],\,\,s\in [l_0,l],
\label{3_14}
\end{equation}
\begin{equation}
\frac{\kappa_0+1}{\mu_0}g'_0(s)=-\frac{\kappa+1}{\mu}g'(s),\,\,s\in [l_0,l].
\label{3_15}
\end{equation}

On the other hand, using the formulas (\ref{3_1}), (\ref{3_2}) and (\ref{3_3}), (\ref{3_4}), it is possible to express the stresses $(\sigma_n+i\tau_n)^{\pm}(t)$ and the derivatives of the displacements $\frac{d}{dt}(u_1+iu_2)^{\pm}(t)$ on the line $L_0\cup L$ through the functions $g'(t)$, $q(t)$:
\begin{equation}
\label{3_16}
(\sigma_n+i\tau_n)^{\pm}(t)=\pm q(t)+\frac{1}{2\pi}\int_{L_0\cup L}\left( \frac{1}{\tau-t}+\frac{1}{\bar{\tau}-\bar{t}}\frac{d\bar{t}}{dt}\right)g'(\tau)d\tau+
\end{equation}
$$
\frac{1}{2\pi}\int_{L_0\cup L}\left(\frac{1}{\bar{\tau}-\bar{t}}-\frac{\tau-t}{(\bar{\tau}-\bar{t})^2}\frac{d\bar{t}}{dt}\right)\overline{g'(\tau)d\tau}+
$$
$$
\frac{(\kappa+1)^{-1}}{\pi i}\int_{L_0\cup L}\left(\frac{1}{\tau-t}-\frac{\kappa}{\bar{\tau}-\bar{t}}\frac{d\bar{t}}{dt} \right)q(\tau)d\tau-
$$
$$
\frac{(\kappa+1)^{-1}}{\pi i}\int_{L_0\cup L}\left(\frac{1}{\bar{\tau}-\bar{t}}-\frac{\tau-t}{(\bar{\tau}-\bar{t})^2}\frac{d\bar{t}}{dt}\right)\overline{q(\tau)d\tau}+2\R \Gamma+\bar{\Gamma}'\frac{d\bar{t}}{dt},\,\,\,t\in L_0\cup L,
$$
$$
2\mu\frac{d}{dt}(u_1+iu_2)^{\pm}(t)=
\pm\frac{i(\kappa+1)}{2}g'(t)+\frac{1}{2\pi}\int_{L_0\cup L}\left( \frac{\kappa}{\tau-t}-\frac{1}{\bar{\tau}-\bar{t}}\frac{d\bar{t}}{dt}\right)g'(\tau)d\tau-
$$
\begin{equation}
\label{3_17}
\frac{1}{2\pi}\int_{L_0\cup L}\left(\frac{1}{\bar{\tau}-\bar{t}}-\frac{\tau-t}{(\bar{\tau}-\bar{t})^2}\frac{d\bar{t}}{dt}\right)\overline{g'(\tau)d\tau}+
\end{equation}
$$
\frac{(\kappa+1)^{-1}}{\pi i}\int_{L_0\cup L}\left(\frac{\kappa}{\tau-t}+\frac{\kappa}{\bar{\tau}-\bar{t}}\frac{d\bar{t}}{dt} \right)q(\tau)d\tau+
$$
$$
\frac{(\kappa+1)^{-1}}{\pi i}\int_{L_0\cup L}\left(\frac{1}{\bar{\tau}-\bar{t}}-\frac{\tau-t}{(\bar{\tau}-\bar{t})^2}\frac{d\bar{t}}{dt}\right)\overline{q(\tau)d\tau}+\kappa \Gamma-\bar{\Gamma}-\bar{\Gamma}'\frac{d\bar{t}}{dt}, \,\,\,t\in L_0\cup L.
$$

Similar formulas for the values $(\sigma_n+i\tau_n)^{\pm}_0(t)$ and $2\mu_0\frac{d}{dt}(u_1+iu_2)^{\pm}_0(t)$ can be obtained by taking all of the parameters in (\ref{3_16}), (\ref{3_17}) with subindex ``$0$" and ommiting the terms containing $\Gamma$ and $\Gamma'$.

Observe, that formally extending the inclusion $S_0$ and the plate $S$ using the formulas (\ref{3_9}) allows us to write the boundary conditions (\ref{2_1}), (\ref{2_5})-(\ref{2_7}) in the form (\ref{3_12a})-(\ref{3_15}), and, hence, to avoid the differentiation of the singular integrals in the formulas (\ref{3_16}), (\ref{3_17}) which would be necessary otherwise. To make sure that the equations (\ref{3_9}) hold, it is sufficient to satisfy the conditions (\ref{3_10a}) which lead to a relatively simple system of singular integral equations:
\begin{equation}
-\frac{i(\kappa_0+1)}{2}g'_0(t)+\frac{1}{2\pi}\int_{L_0\cup L}\left(\frac{\kappa_0-1}{\tau-t}-k_1(t,\tau) \right)g'_0(\tau)d\tau-
\label{3_23}
\end{equation}
$$
\frac{1}{2\pi}\int_{L_0\cup L}k_2(t,\tau)\overline{g'_0(\tau)d\tau}+
\frac{\kappa_0(\kappa_0+1)^{-1}}{\pi i}\int_{L_0\cup L}\left(\frac{2}{\tau-t}+k_1(t,\tau) \right)q_0(\tau)d\tau+
$$
$$
\frac{(\kappa_0+1)^{-1}}{\pi i}\int_{L_0\cup L} k_2(t,\tau)\overline{q_0(\tau)d\tau}
=0,\,\,\,t\in L_0\cup L,
$$
\begin{equation}
\frac{i(\kappa+1)}{2}g'(t)+\frac{1}{2\pi}\int_{L_0\cup L}\left(\frac{\kappa-1}{\tau-t}-k_1(t,\tau) \right)g'(\tau)d\tau-
\label{3_24}
\end{equation}
$$
\frac{1}{2\pi}\int_{L_0\cup L}k_2(t,\tau)\overline{g'(\tau)d\tau}+\frac{\kappa(\kappa+1)^{-1}}{\pi i}\int_{L_0\cup L}\left(\frac{2}{\tau-t}+k_1(t,\tau) \right)q(\tau)d\tau+
$$
$$
\frac{(\kappa+1)^{-1}}{\pi i}\int_{L_0\cup L} k_2(t,\tau)\overline{q(\tau)d\tau}
+\kappa\Gamma-\bar{\Gamma}-\bar{\Gamma}'\frac{d\bar{t}}{dt}=0,\,\,\,t\in L_0\cup L,
$$
where $k_1(t,\tau)$ and $k_2(t,\tau)$ are the regular kernels defined by the formulas:
$$
k_1(t,\tau)=\frac{d}{dt}\ln \left(\frac{\tau-t}{\bar{\tau}-\bar{t}} \right)=-\frac{1}{\tau-t}+\frac{1}{\bar{\tau}-\bar{t}}\frac{d\bar{t}}{dt},
$$
$$
k_2(t,\tau)=-\frac{d}{dt}\frac{\tau-t}{\bar{\tau}-\bar{t}}=\frac{1}{\bar{\tau}-\bar{t}}-\frac{\tau-t}{(\bar{\tau}-\bar{t})^2}\frac{d\bar{t}}{dt}.
$$

Observe, that integrating the right-hand side of the equation (\ref{3_24}) along the contour $L_0\cup L$ produces an identical zero. Thus, the general solution of the system (\ref{3_23}), (\ref{3_24}) contains two real constants which can be fixed by stating two additional real conditions. It is required from a physical viewpoint that the displacements are single-valued along the contour $L_0$ which means that the relative displacements of the crack tips should be the same traced along both left- and right-hand-side banks of the crack. This condition is not satisfied automatically and needs to be stated additionally as a part of the solution:
\begin{equation}
\int_{L_0}\frac{d}{dt}(u_1+iu_2)^+_0(t)dt=\int_{L_0}\frac{d}{dt}(u_1+iu_2)^-(t)dt.
\label{3_25}
\end{equation}
The equation (\ref{3_25}) can be easily written in terms of the unknown functions using the conditions (\ref{3_10a}):
\begin{equation}
\int_{L_0}\left(\frac{\kappa_0+1}{\mu_0}g'_0(\tau)+\frac{\kappa+1}{\mu}g'(\tau) \right)d\tau=0.
\label{3_26}
\end{equation}
The condition (\ref{3_26}) can be incorporated into the equation (\ref{3_24}) by adding an extra term of the form $\int_{L_0}\left(\frac{\kappa_0+1}{\mu_0}g'_0(\tau)+\frac{\kappa+1}{\mu}g'(\tau) \right)d\tau\frac{|dt|}{dt}$ to the left-hand side of the equation. Then the equation (\ref{3_24}) becomes:
\begin{equation}
\frac{i(\kappa+1)}{2}g'(t)+\frac{1}{2\pi}\int_{L_0\cup L}\left(\frac{\kappa-1}{\tau-t}-k_1(t,\tau) \right)g'(\tau)d\tau-
\label{3_27}
\end{equation}
$$
\frac{1}{2\pi}\int_{L_0\cup L}k_2(t,\tau)\overline{g'(\tau)d\tau}+\frac{\kappa(\kappa+1)^{-1}}{\pi i}\int_{L_0\cup L}\left(\frac{2}{\tau-t}+k_1(t,\tau) \right)q(\tau)d\tau+
$$
$$
\frac{(\kappa+1)^{-1}}{\pi i}\int_{L_0\cup L} k_2(t,\tau)\overline{q(\tau)d\tau}
+\int_{L_0}\left(\frac{\kappa_0+1}{\mu_0}g'_0(\tau)+\frac{\kappa+1}{\mu}g'(\tau) \right)d\tau\frac{|dt|}{dt}+
$$
$$
\kappa\Gamma-\bar{\Gamma}-\bar{\Gamma}'\frac{d\bar{t}}{dt}=0,\,\,\,t\in L_0\cup L.
$$

\setcounter{equation}{0}
\section{System of singular integro-differential equations}

Observe, that the solutions $g'_0(t)$, $q_0(t)$, $g'(t)$, $q(t)$ to the system (\ref{3_23}), (\ref{3_27}) are sought in the class of H\"older continuous functions on the curves $L_0$, $L$ and can have at most integrable power singularities at the end points of the lines $L_0$, $L$. Hence, since the contour $L_0\cup L$ is closed it is possible to use the inversion formula for a singular integral of Cauchy type \cite{MikhPros1986}:
\begin{equation}
S^2_{L_0\cup L}=I,
\label{4_1}
\end{equation}
where $S_{L_0\cup L}$ denotes the Cauchy singular integral operator on the curve $L_0\cup L$, and $I$ is the identity operator.

Separating the regular and the singular parts in the equations (\ref{3_23}), (\ref{3_27}) and applying the formula (\ref{4_1}), obtain:
\begin{equation}
\frac{2\kappa_0}{\kappa_0+1}q_0(t)+\frac{i(\kappa_0-1)}{2}g'_0(t)-\frac{\kappa_0+1}{2\pi}\int_{L_0\cup L}\frac{g'_0(\tau)d\tau}{\tau-t}=
\label{4_2}
\end{equation}
$$
\frac{\kappa_0(\kappa_0+1)^{-1}}{\pi^2}\int_{L_0\cup L}q_0(\tau)d\tau\int_{L_0\cup L}\frac{k_1(\tau_1,\tau)d\tau_1}{\tau_1-t}+
$$
$$
\frac{(\kappa_0+1)^{-1}}{\pi^2}\int_{L_0\cup L}\overline{q_0(\tau)d\tau}\int_{L_0\cup L}\frac{k_2(\tau_1,\tau)d\tau_1}{\tau_1-t}+
$$
$$
\frac{1}{2\pi^2i}\int_{L_0\cup L}g'_0(\tau)d\tau\int_{L_0\cup L}\frac{k_1(\tau_1,\tau)d\tau_1}{\tau_1-t}-
$$
$$
\frac{1}{2\pi^2i}\int_{L_0\cup L}\overline{g'_0(\tau)d\tau}\int_{L_0\cup L}\frac{k_2(\tau_1,\tau)d\tau_1}{\tau_1-t},\,\,\,t\in L_0\cup L,
$$
\begin{equation}
\frac{2\kappa}{\kappa+1}q(t)+\frac{i(\kappa-1)}{2}g'(t)+\frac{\kappa+1}{2\pi}\int_{L_0\cup L}\frac{g'(\tau)d\tau}{\tau-t}=
\label{4_3}
\end{equation}
$$
\frac{\kappa(\kappa+1)^{-1}}{\pi^2}\int_{L_0\cup L}q(\tau)d\tau\int_{L_0\cup L}\frac{k_1(\tau_1,\tau)d\tau_1}{\tau_1-t}+
$$
$$
\frac{(\kappa+1)^{-1}}{\pi^2}\int_{L_0\cup L}\overline{q(\tau)d\tau}\int_{L_0\cup L}\frac{k_2(\tau_1,\tau)d\tau_1}{\tau_1-t}+
$$
$$
\frac{1}{2\pi^2i}\int_{L_0\cup L}g'(\tau)d\tau\int_{L_0\cup L}\frac{k_1(\tau_1,\tau)d\tau_1}{\tau_1-t}-\frac{1}{2\pi^2i}\int_{L_0\cup L}\overline{g'(\tau)d\tau}\int_{L_0\cup L}\frac{k_2(\tau_1,\tau)d\tau_1}{\tau_1-t}+
$$
$$
\frac{1}{\pi i}\int_{L_0}\left(\frac{\kappa_0+1}{\mu_0}g'_0(\tau)+\frac{\kappa+1}{\mu}g'(\tau) \right)d\tau\int_{L_0\cup L}\frac{|d\tau|}{\tau-t}-
$$
$$
\kappa\Gamma+\bar{\Gamma}+\frac{\bar{\Gamma}'}{\pi i}\int_{L_0\cup L}\frac{d\bar{\tau}}{\tau-t},\,\,\,t\in L_0\cup L.
$$

By separating the real and the imaginary parts of the equations (\ref{3_23}), (\ref{3_27}), (\ref{4_2}), (\ref{4_3}), eliminating the term $\int_{L_0\cup L}\frac{\I g'_0(\tau)d\tau}{\tau-t}$, and also transforming the equations into the natural parameter $s$, obtain:
$$
\frac{2\kappa_0}{\kappa_0+1}\R q_0(s_0)-\frac{\kappa_0-1}{2}\I g'_0(s_0)-\frac{\kappa_0+1}{2\pi}\int_0^l\frac{\R g'_0(s)ds}{s-s_0}=
$$
\begin{equation}
M_1(g'_0,q_0)(s_0),\,\,\,s_0\in [0,l],
\label{4_4}
\end{equation}
\begin{equation}
\frac{\kappa_0-1}{\kappa_0+1}\I q_0(s_0)-\R g'_0(s_0)-\frac{1}{\pi}\int_{0}^l\frac{\R q_0(s)ds}{s-s_0}=M_2(g'_0,q_0)(s_0),\,\,\,s_0\in [0,l],
\label{4_5}
\end{equation}
\begin{equation}
\frac{2\kappa}{\kappa+1}\R q(s_0)-\frac{\kappa-1}{2}\I g'(s_0)+\frac{\kappa+1}{2\pi}\int_0^l\frac{\R g'(s)ds}{s-s_0}=
\label{4_6}
\end{equation}
$$
M_3(g',q)(s_0),\,\,\,s_0\in [0,l],
$$
\begin{equation}
\frac{\kappa-1}{\kappa+1}\I q(s_0)-\R g'(s_0)+\frac{1}{\pi}\int_{0}^l\frac{\R q(s)ds}{s-s_0}=M_4(g',q)(s_0),\,\,\,s_0\in [0,l],
\label{4_7}
\end{equation}
where $M_1(g'_0,q_0)(s_0)$, $M_2(g'_0,q_0)(s_0)$, $M_3(g',q)(s_0)$, $M_4(g',q)(s_0)$ are the real-valued regular terms which can be easily obtained from the equations (\ref{4_2}), (\ref{4_3}).

The equations (\ref{4_4})-(\ref{4_7}) are valid on both of the segments $[0,l_0]$, $[l_0,l]$ and contain eight real unknown functions $\R g'_0(s_0)$, $\I g'_0(s_0)$, $\R q_0(s_0)$, $\I q_0(s_0)$, $\R g'(s_0)$, $\I g'(s_0)$, $\R q(s_0)$, $\I q(s_0)$ which are defined on both segments $[0,l_0]$, $[l_0,l]$, and satisfy the functional relationships (\ref{3_12a})-(\ref{3_15}).

Introduce the new real unknown functions by the following formulas:
\begin{equation}
\chi_1(s_0)=\frac{\gamma^+}{4\mu_0}(\kappa_0+1)\frac{d}{ds_0}\left[\varrho(s_0)\I g'_0(s_0)+\R g''_0(s_0) \right],\,\,\,s_0\in [0,l_0],
\label{4_8}
\end{equation}
\begin{equation}
\chi_2(s_0)=\I g'_0(s_0),\,\,\,s_0\in [0,l_0],
\label{4_9}
\end{equation}
\begin{equation}
\chi_3(s_0)=\frac{\gamma^-}{4\mu}(\kappa+1)\frac{d}{ds_0}\left[\varrho(s_0)\I g'(s_0)+\R g''(s_0) \right],\,\,\,s_0\in [0,l_0],
\label{4_10}
\end{equation}
\begin{equation}
\chi_4(s_0)=\I g'(s_0),\,\,\,s_0\in [0,l_0],
\label{4_11}
\end{equation}
\begin{equation}
\chi_5(s_0)=\I q_0(s_0)+\I q(s_0),\,\,\,s_0\in [l_0,l],
\label{4_12}
\end{equation}
\begin{equation}
\chi_6(s_0)=\I g'_0(s_0),\,\,\,s_0\in [l_0,l],
\label{4_13}
\end{equation}
\begin{equation}
\chi_7(s_0)=\frac{\kappa_0}{\mu_0(\kappa_0+1)}\R q_0(s_0)+\frac{\kappa}{\mu(\kappa+1)}\R q(s_0),\,\,\,s_0\in [l_0,l],
\label{4_14}
\end{equation}
\begin{equation}
\chi_8(s_0)=\frac{\kappa_0-1}{\kappa_0+1}\I q_0(s_0)-\frac{\kappa-1}{\kappa+1}\I q(s_0),\,\,\,s_0\in [l_0,l].
\label{4_15}
\end{equation}

Substituting the identities (\ref{4_8})-(\ref{4_15}) into the equations (\ref{3_12a})-(\ref{3_15}) allows us to express the unknown functions $\R g'_0(s_0)$, $\I g'_0(s_0)$, $\R q_0(s_0)$, $\I q_0(s_0)$, $\R g'(s_0)$, $\I g'(s_0)$, $\R q(s_0)$, $\I q(s_0)$ through the new unknown functions $\chi_j(s_0)$, $j=1,2,\ldots,8$:
\begin{equation}
\I q_0(s_0)=\chi_1(s_0)+\frac{1}{2}\I f_1(s_0),\,\,\,s_0\in [0,l_0],
\label{4_16}
\end{equation}
\begin{equation}
\I g'_0(s_0)=\chi_2(s_0),\,\,\,s_0\in [0,l_0],
\label{4_17}
\end{equation}
\begin{equation}
\R q_0(s_0)=\varrho(s_0)\left(\int_0^{s_0}\chi_1(s)ds+C_1\right)+\frac12\R f_1(s_0),\,\,\,s_0\in [0,l_0],
\label{4_18}
\end{equation}
\begin{equation}
\R g'_0(s_0)=\frac{4\mu_0}{\gamma^+(\kappa_0+1)}\left(\int_0^{s_0}(s_0-s)\chi_1(s)ds+C_1s_0\right)-
\label{4_19}
\end{equation}
$$
\int_0^{s_0}\varrho(s)\chi_2(s)ds+C_2,\,\,\,s_0\in [0,l_0],
$$
\begin{equation}
\I q(s_0)=\chi_3(s_0)-\frac{1}{2}\I f_2(s_0),\,\,\,s_0\in [0,l_0],
\label{4_20}
\end{equation}
\begin{equation}
\I g'(s_0)=\chi_4(s_0),\,\,\,s_0\in [0,l_0],
\label{4_21}
\end{equation}
\begin{equation}
\R q(s_0)=\varrho(s_0)\left(\int_0^{s_0}\chi_3(s)ds+C_3\right)-\frac12\R f_2(s_0),\,\,\,s_0\in [0,l_0],
\label{4_22}
\end{equation}
\begin{equation}
\R g'(s_0)=\frac{4\mu}{\gamma^-(\kappa+1)}\left(\int_0^{s_0}(s_0-s)\chi_3(s)ds+C_3s_0\right)-
\label{4_23}
\end{equation}
$$
\int_0^{s_0}\varrho(s)\chi_4(s)ds+C_4,\,\,\,s_0\in [0,l_0],
$$
\begin{equation}
\I q_0(s_0)=\left(\frac{\kappa_0-1}{\kappa_0+1}+\frac{\kappa-1}{\kappa+1}\right)^{-1}\left(\frac{\kappa-1}{\kappa+1}\chi_5(s_0)+\chi_8(s_0)\right),\,\,\,s_0\in [l_0,l],
\label{4_24}
\end{equation}
\begin{equation}
\I q(s_0)=\left(\frac{\kappa_0-1}{\kappa_0+1}+\frac{\kappa-1}{\kappa+1}\right)^{-1}\left(\frac{\kappa_0-1}{\kappa_0+1}\chi_5(s_0)-\chi_8(s_0)\right),\,\,\,s_0\in [l_0,l],
\label{4_25}
\end{equation}
$$
\R q_0(s_0)=\left(\frac{\kappa_0}{\mu_0(\kappa_0+1)}-\frac{\kappa}{\mu(\kappa+1)}\right)^{-1}\left[-\frac{\kappa\varrho(s_0)}{\mu(\kappa+1)}\left(\int_0^{s_0}\chi_5(s)ds+C_5\right)+\right.
$$
\begin{equation}
\left.\chi_7(s_0)\right],\,\,\,s_0\in [l_0,l],
\label{4_26}
\end{equation}
$$
\R q(s_0)=\left(\frac{\kappa_0}{\mu_0(\kappa_0+1)}-\frac{\kappa}{\mu(\kappa+1)}\right)^{-1}\left[\frac{\kappa_0\varrho(s_0)}{\mu_0(\kappa_0+1)}\left(\int_0^{s_0}\chi_5(s)ds+C_5\right)-\right.
$$
\begin{equation}
\left.\chi_7(s_0)\right],\,\,\,s_0\in [l_0,l],
\label{4_27}
\end{equation}
\begin{equation}
\I g'_0(s_0)=\chi_6(s_0),\,\,\,s_0\in [l_0,l],
\label{4_28}
\end{equation}
\begin{equation}
\I g'(s_0)=-\frac{\mu(\kappa_0+1)}{\mu_0(\kappa+1)}\chi_6(s_0),\,\,\,s_0\in [l_0,l],
\label{4_29}
\end{equation}
\begin{equation}
\R g'_0(s_0)=\frac{4\mu_0}{\gamma^i(\kappa_0+1)}\left(\int_0^{s_0}(s_0-s)\chi_5(s)ds+C_5s_0\right)-
\label{4_30}
\end{equation}
$$
\int_0^{s_0}\varrho(s)\chi_6(s)ds+C_6,\,\,\,s_0\in [l_0,l],
$$
\begin{equation}
\R g'(s_0)=-\frac{\mu(\kappa_0+1)}{\mu_0(\kappa+1)}\left[\frac{4\mu_0}{\gamma^i(\kappa_0+1)}\left(\int_0^{s_0}(s_0-s)\chi_5(s)ds+C_5s_0\right)\right.-
\label{4_31}
\end{equation}
$$
\left.\int_0^{s_0}\varrho(s)\chi_6(s)ds+C_6\right],\,\,\,s_0\in [l_0,l].
$$
Observe that the formulas (\ref{4_16})-(\ref{4_31}) contain six real unknown parameters $C_1$, $C_2$, ..., $C_6$ which will be fixed later using additional conditions. 

It is possible to regularize the system of the singular integro-differential equations (\ref{4_4})-(\ref{4_7}) in terms of the new unknowns $\chi_j(s)$, $j=1,2,\ldots,8$, defined by the formulas (\ref{4_8})-(\ref{4_15}). Substituting the unknowns (\ref{4_8})-(\ref{4_15}) into the system (\ref{4_4})-(\ref{4_7}), solving the system with respect to $\I g'_0(s_0)$, $\I q_0(s_0)$, $\I g'(s_0)$, $\I q(s_0)$, and changing the order of integration where appropriate, reduces the system of the singular integro-differential equations (\ref{4_4})-(\ref{4_7}) on the interval $[0, l_0]$ to the system of weakly-singular integral equations solved with respect to the unknowns $\chi_j(s)$, $j=1,2,3,4$.

A slightly different approach is necessary to regularize the system (\ref{4_4})-(\ref{4_7}) on the interval $[l_0,l]$. First, multiply the equation (\ref{4_4}) by $(\kappa_0+1)/(2\kappa_0)$, multiply the equation (\ref{4_6}) by $(\kappa+1)/(2\kappa)$, and add:
$$
\R q_0(s_0)+\R q(s_0)-\frac{\kappa_0^2-1}{2\kappa_0}\I g'_0(s_0)-\frac{\kappa^2-1}{4\kappa}\I g'(s_0)-
$$
\begin{equation}
\frac{(\kappa_0+1)^2}{4\pi\kappa_0}\int_0^l\frac{\R g'_0(s)ds}{s-s_0}+\frac{(\kappa+1)^2}{4\pi\kappa}\int_0^l\frac{\R g'(s)ds}{s-s_0}=
\label{4_32}
\end{equation}
$$
\frac{\kappa_0+1}{4\kappa_0}M_1(g'_0,q_0)(s_0)+\frac{\kappa+1}{4\kappa}M_3(g',q)(s_0),\,\,\,s\in [l_0,l].
$$
Next, multiply the equation (\ref{4_4}) by $(\kappa_0+1)/(2\mu_0(\kappa_0-1))$, multiply the equation (\ref{4_6}) by $(\kappa+1)/(2\mu(\kappa-1))$, and add:
\begin{equation}
\frac{\kappa_0}{\mu_0(\kappa_0-1)}\R q_0(s_0)+\frac{\kappa}{\mu(\kappa-1)}\R q(s_0)-\frac{(\kappa_0+1)^2}{4\pi\mu_0(\kappa_0-1)}\int_0^l\frac{\R g'_0(s)ds}{s-s_0}+
\label{4_33}
\end{equation}
$$
\frac{(\kappa+1)^2}{4\pi\mu(\kappa-1)}\int_0^L\frac{\R g'(s)ds}{s-s_0}=\frac{\kappa_0+1}{2\mu_0(\kappa_0-1)}M_1(g'_0,q_0)(s_0)+\frac{\kappa+1}{2\mu(\kappa-1)}M_3(g',q)(s_0),
$$
$$
s_0\in [l_0,l].
$$
To obtain the third equation of the new system, subtract the equation (\ref{4_7}) from the equation (\ref{4_5}):
$$
\frac{\kappa_0-1}{\kappa_0+1}\I q_0(s_0)-\frac{\kappa-1}{\kappa+1}\I q(s_0)-\R g'_0(s_0)+\R g'(s_0)-
$$
\begin{equation}
\frac{1}{\pi}\int_0^l\frac{(\R q_0(s)+\R q(s))ds}{s-s_0}=M_2(g'_0,q_0)-M_4(g',q)(s_0),\,\,\,s_0\in [l_0,l].
\label{4_34}
\end{equation}
Finally, multiply the equation (\ref{4_5}) by $(\kappa_0+1)/(\kappa_0-1)$, multiply the equation (\ref{4_7}) by $(\kappa+1)/(\kappa-1)$, and add, to obtain the last equation of the modified system:
$$
\I q_0(s_0)+\I q(s_0)-\frac{\kappa_0+1}{\kappa_0-1}\R g'_0(s_0)-\frac{\kappa+1}{\kappa-1}\R g'(s_0)-
$$
\begin{equation}
\frac{\kappa_0+1}{\pi(\kappa_0-1)}\int_0^l\frac{\R q_0(s)ds}{s-s_0}+\frac{\kappa+1}{\pi(\kappa-1)}\int_0^l\frac{\R q(s)ds}{s-s_0}=
\label{4_35}
\end{equation}
$$
\frac{\kappa_0+1}{\kappa_0-1}M_2(g'_0,q_0)+\frac{\kappa+1}{\kappa-1}M_4(g',q)(s_0),\,\,\,s_0\in [l_0,l].
$$

Notice that the equations (\ref{4_32})-(\ref{4_34}), after substituting the unknowns (\ref{4_12})-(\ref{4_15}) and making the appropriate changes of the order of integration, produce weakly-singular integral equations solved with respect to the unknown functions $\chi_6(s_0)$, $\chi_7(s_0)$, $\chi_8(s_0)$. Finally, substituting the equations (\ref{4_32}), (\ref{4_33}) into the equation (\ref{4_35}) gives us one more weakly-singular equation solved with respect to the unknown function $\chi_5(s_0)$. The procedure described above works for all the material parameters such as $\mu_0\kappa(\kappa_0+1)\neq \mu\kappa_0(\kappa+1)$. In the special case $\mu_0\kappa(\kappa_0+1)=\mu\kappa_0(\kappa+1)$ the regularization of the system (\ref{4_4})-(\ref{4_7}) still can be accomplished by choosing different unknowns (\ref{4_12})-(\ref{4_15}).

Hence, the system of the singular integro-differential equations (\ref{4_4})-(\ref{4_7}) on the intervals $[0,l_0]$ and $[l_0,l]$ can be reduced to the system of eight weakly-singular integral equations of the second kind with respect to the unknowns $\chi_j(s)$, $j=1,2,\ldots,8$, defined by the formulas (\ref{4_8})-(\ref{4_15}). The theory of these systems is similar to the theory of the Fredholm equations of the second kind. Thus, the system has a unique solution for almost all values of the material and geometric parameters.

\setcounter{equation}{0}
\section{Free constants, additional conditions and singularities at the crack tips}

The formulas (\ref{4_16})-(\ref{4_31}) contain six real unknown constants $C_1$, $C_2$, ..., $C_6$ which need to be fixed by specifying six real additional conditions. Two of these conditions can be obtained by requiring that the total force applied to the inclusion and the matrix from both sides of the line $L_0\cup L$ is equal to zero:
$$
\int_{L_0\cup L}\left((\sigma_n+i\tau_n)^+_0(t)+(\sigma_n+i\tau_n)^-(t)\right)dt=0.
$$
Using the equations (\ref{3_5}), (\ref{3_7}), (\ref{3_11}), the last condition can be rewritten in the form:
\begin{equation}
\int_{L_0\cup L}(q_0(t)-q(t))dt=0.
\label{5_1}
\end{equation}

There are several ways to fix the remaining four real conditions. Here, it is assumed that
\begin{equation}
\R g'_0(0+)=\R g'_0(l-0),\,\,\,\R g'_0(l_0-0)=\R g'_0(l_0+0),
\label{5_2}
\end{equation}
\begin{equation}
\R g'(0+)=\R g'(l-0),\,\,\,\R g'(l_0-0)=\R g'(l_0+0).
\label{5_3}
\end{equation}
As it will be seen below, the imposed conditions (\ref{5_2}), (\ref{5_3}) guarantee the solution to the mechanical problem which allows for at most logarithmic singularities of the stresses and the derivatives of the displacements at the crack tips.

Consider the integrals
\begin{equation}
\int_{L_0\cup L}\frac{\R g'_0(\tau)d\tau}{\tau-t},\,\,\,\int_{L_0\cup L}\frac{\R g'(\tau)d\tau}{\tau-t}. 
\label{5_4}
\end{equation}
Due to the conditions (\ref{5_2}), (\ref{5_3}) it follows that the functions $\R g'_0(t)$, $\R g'(t)$ are continuous on the closed contour $L_0\cup L$ including the end-points of the curves $L_0$ and $L$. Hence, the integrals (\ref{5_4}) are bounded everywhere on the contour $L_0\cup L$. Recall that according to the formulas (\ref{4_18}), (\ref{4_22}) the limits $\R q_0(0+)$, $\R q_0(l_0-0)$, $\R q(0+)$, $\R q(l_0-0)$ exist and are finite. Then, by using the equations (\ref{4_4}), (\ref{4_6}) on the interval $[0,l_0]$, it follows that the limits $\I g'_0(0+)$, $\I g'_0(l_0-0)$, $\I g'(0+)$, $\I g'(l_0-0)$ exist and are finite as well. Also, from the equation (\ref{4_33}) it follows that $\frac{\kappa_0}{\mu_0(\kappa_0-1)}\R q_0(s_0)+\frac{\kappa}{\mu(\kappa-1)}\R q(s_0)$ is bounded as $s_0\to l_0+0$ and $s_0 \to l-0$. Combined with the fact that the expression $\R q_0(s_0)+\R q(s_0)$ is bounded at these points as well according to the formulas (\ref{4_26}), (\ref{4_27}), it follows that the limits $\R q_0(l_0+0)$, $\R q_0(l-0)$, $\R q(l_0+0)$, $\R q(l-0)$ exist and are finite. Then from the equations (\ref{4_4}), (\ref{4_6}) on the interval $[l_0,l]$, it follows that the values $\I g'_0(l_0+0)$, $\I g'_0(l-0)$, $\I g'(l_0+0)$, $\I g'(l-0)$ are finite as well. Finally, since the conditions 
$\R q_0(0+)=\R q_0(l-0)$, $\R q_0(l_0-0)=\R q_0(l_0+0)$, $\R q(0+)=\R q(l-0)$, $\R q(l_0-0)=\R q(l_0+0)$ are not necessary satisfied, the equations (\ref{4_5}), (\ref{4_7}) suggest that the functions $\I q_0(s_0)$, $\I q(s_0)$ may still possess at most logarithmic singularities at $s_0=0$, $s_0=l_0$ and $s_0=l$.

The discussed behavior of the functions $g'_0(s)$, $q_0(s)$, $g'(s)$, $q(s)$ together with the conditions (\ref{3_5})-(\ref{3_8}) and (\ref{3_10}), (\ref{3_11}) guarantees that the derivatives of the displacements ${d(u_1+iu_2)}/{dt}$ are bounded at the crack tips, the normal stresses $\sigma_n$ are also bounded, while the shear stresses $\tau_n$ may have at most logarithmic singularities at the crack tips. This is a significant improvement over the LEFM theory since both power singularities of the order $1/2$ and oscillating singularities are eliminated by the introduction of the curvature-dependent surface tension on the boundary of the crack. This conclusion is in agreement with the results of the papers \cite{SendovaWalton2010}, \cite{Zemlyanova2013}, \cite{Zemlyanova2012}.

\setcounter{equation}{0}
\section{Numerical solution of the system of the singular integro-differential equations}

The system of the weakly-singular integral equations, obtained by the regularization of the system (\ref{3_23}), (\ref{3_27}) with the conditions (\ref{3_12a})-(\ref{3_15}), in general, has a very complicated form and is not convenient for using in numerical computations. Therefore, it is customary to solve numerically the initial system of the singular integro-differential equations (\ref{3_23}), (\ref{3_27}) together with the conditions (\ref{5_1})-(\ref{5_3}). There are multiple ways to solve the systems of this type numerically, such as spline collocation methods and representations of unknowns with different special functions \cite{Samoilova2001}, \cite{Schmidt1987}.  In this work we will follow the approach adopted in \cite{Zemlyanova2012}.
The unknown functions will be sought in the form of truncated Taylor series (Taylor polynomials) separately on the segments $[0,l_0]$ and $[l_0,l]$:

\begin{figure}[ht]
	\centering
		\scalebox{0.5}{\includegraphics{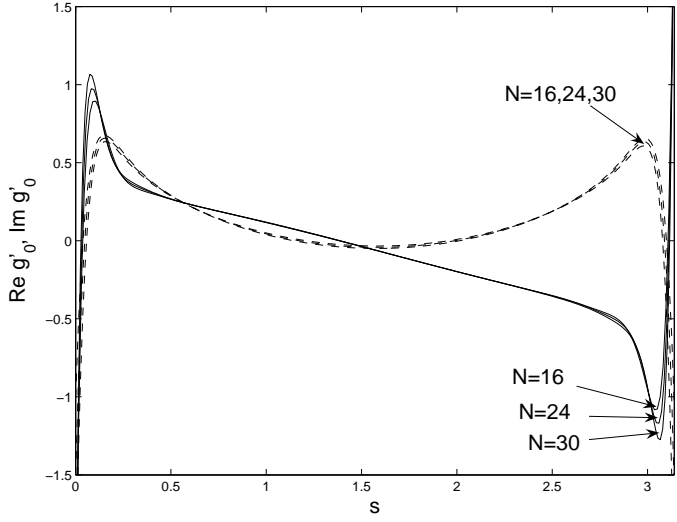}}
			\caption{Graphs of the functions $\R g'_0(s)$ and $\I g'_0(s)$ on the boundary of a semicircular interface fracture $L_0$ for different values of the parameter $N$}
	\label{fig1}
\end{figure}

$$
q_0(s)=\sum_{k=0}^{N+1}a^1_k(s-l_0/2)^k+i\sum_{k=0}^{N}b^1_k(s-l_0/2)^k,\,\,\,s\in [0,l_0],
$$
$$
g'_0(s)=\sum_{k=0}^{N+1}a^2_k(s-l_0/2)^k+i\sum_{k=0}^{N}b^2_k(s-l_0/2)^k,\,\,\,s\in [0,l_0],
$$
$$
q(s)=\sum_{k=0}^{N+1}a^3_k(s-l_0/2)^k+i\sum_{k=0}^{N}b^3_k(s-l_0/2)^k,\,\,\,s\in [0,l_0],
$$
\begin{equation}
g'(s)=\sum_{k=0}^{N+1}a^4_k(s-l_0/2)^k+i\sum_{k=0}^{N}b^4_k(s-l_0/2)^k,\,\,\,s\in [0,l_0],
\label{6_1}
\end{equation}
$$
q_0(s)=\sum_{k=0}^{N+1}a^5_k(s-l_1/2)^k+i\sum_{k=0}^{N}b^5_k(s-l_1/2)^k,\,\,\,s\in [l_0,l],
$$
$$
g'_0(s)=\sum_{k=0}^{N+1}a^6_k(s-l_1/2)^k+i\sum_{k=0}^{N}b^6_k(s-l_1/2)^k,\,\,\,s\in [l_0,l],
$$
$$
q(s)=\sum_{k=0}^{N}a^7_k(s-l_1/2)^k+i\sum_{k=0}^{N}b^7_k(s-l_1/2)^k,\,\,\,s\in [l_0,l],
$$
$$
g'(s)=\sum_{k=0}^{N+1}a^8_k(s-l_1/2)^k+i\sum_{k=0}^{N}b^8_k(s-l_1/2)^k,\,\,\,s\in [l_0,l],
$$
where $a^j_k$ and $b^j_k$ are unknown real coefficients, and $l_1=(l_0+l)/2$. The number of the terms in the formulas (\ref{6_1}) is chosen so that the representations (\ref{6_1}) contain six additional real constants which correspond to the constants of integration $C_1$, $C_2$, .., $C_6$ in the formulas (\ref{4_16})-(\ref{4_31}).

Observe, that from the conditions (\ref{3_15}) it follows that
\begin{equation}
\frac{\kappa_0+1}{\mu_0}a^6_k=-\frac{\kappa+1}{\mu}a^8_k,\, k=1,\ldots,N+1,
\label{6_2}
\end{equation}
$$
\frac{\kappa_0+1}{\mu_0}b^6_k=-\frac{\kappa+1}{\mu}b^8_k,\, k=1,\ldots,N. 
$$
Substituting the formulas (\ref{6_1}) into the equations (\ref{3_23}), (\ref{3_27}), (\ref{3_12a})-(\ref{3_14}) and additional conditions (\ref{5_1})-(\ref{5_3}),  produces a system of equations with respect to the unknowns $a^j_k$, $b^j_k$.  Satisfying the equations of this system only in $N+1$ discrete equally spaced points on the intervals $[0,l_0]$ and $[l_0,l]$, and taking into account the equations (\ref{6_2}), reduces the problem to the solution of a system of linear algebraic equations with $16N+23$ unknowns and equations. This can be done numerically using Matlab software.

\begin{figure}[ht]
	\centering
		\scalebox{0.3}{\includegraphics{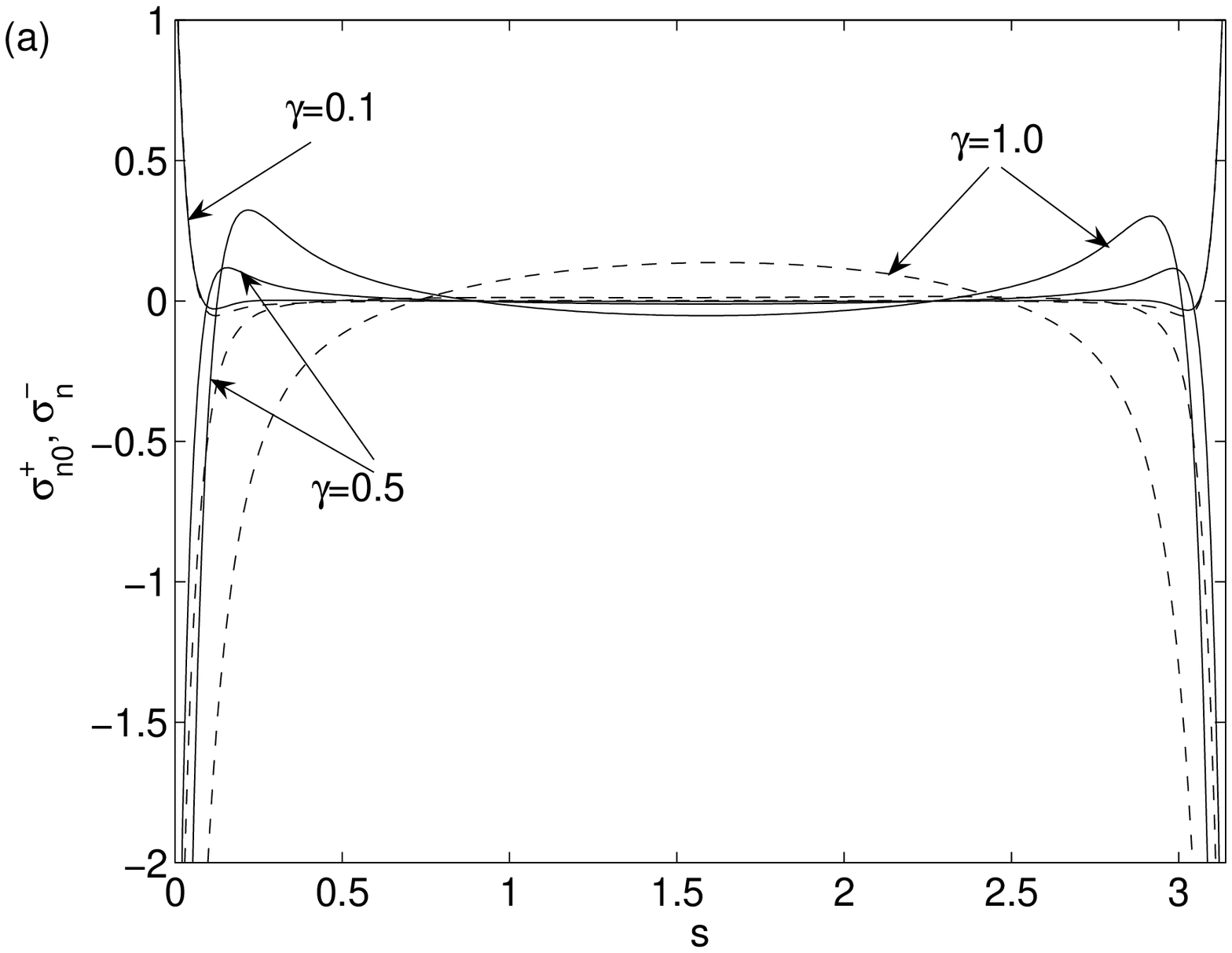} \hspace{1.0cm} \includegraphics{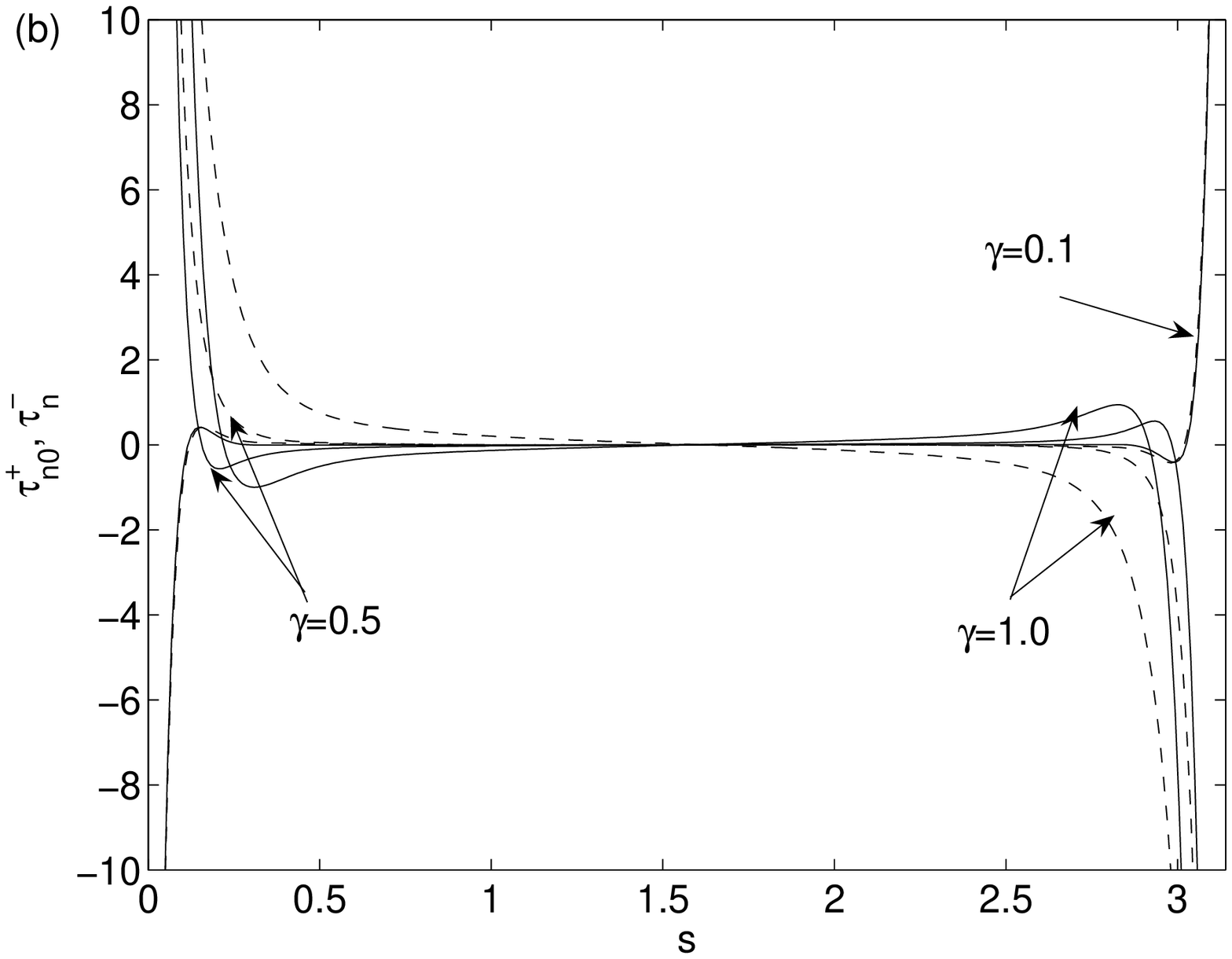}}
		\scalebox{0.3}{\includegraphics{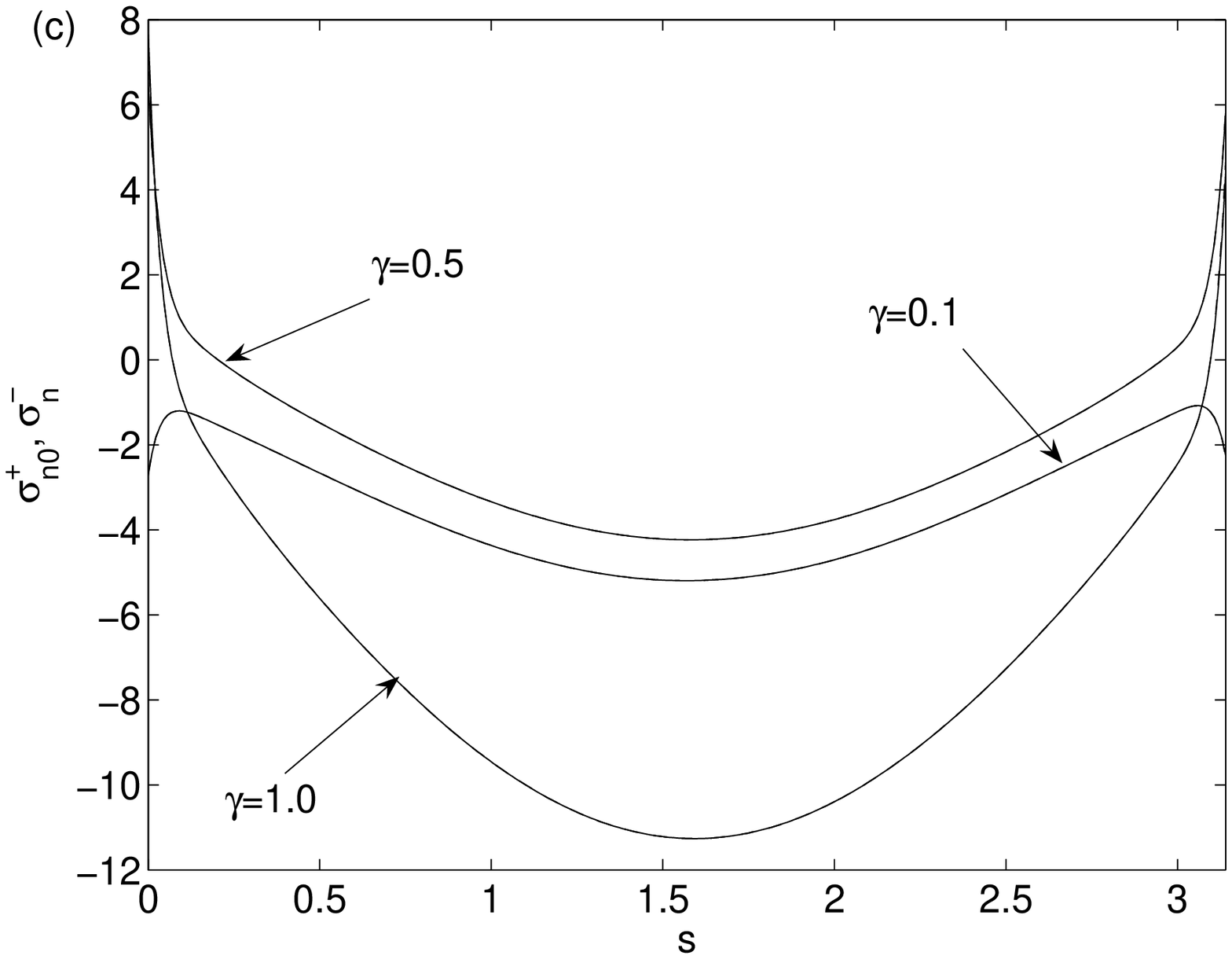} \hspace{1.0cm} \includegraphics{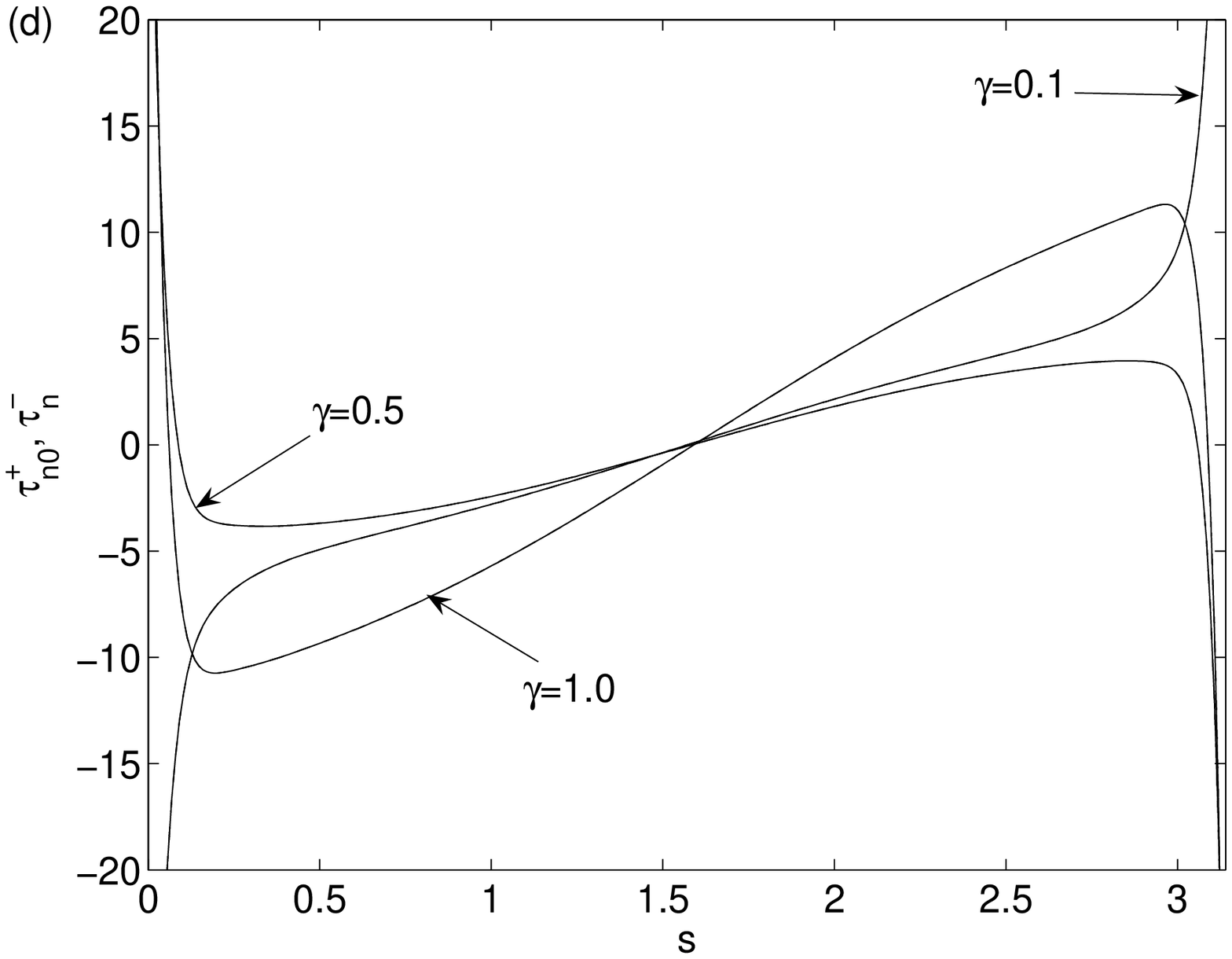}}
			\caption{Graphs of the normal and shear stresses $\sigma^{+}_{n0}$, $\sigma_n^-$, $\tau^{+}_{n0}$, $\tau_n^-$ on the curves $L_0$ and $L$}
	\label{fig2}
\end{figure}

\begin{figure}[ht]
	\centering
		\scalebox{0.3}{\includegraphics{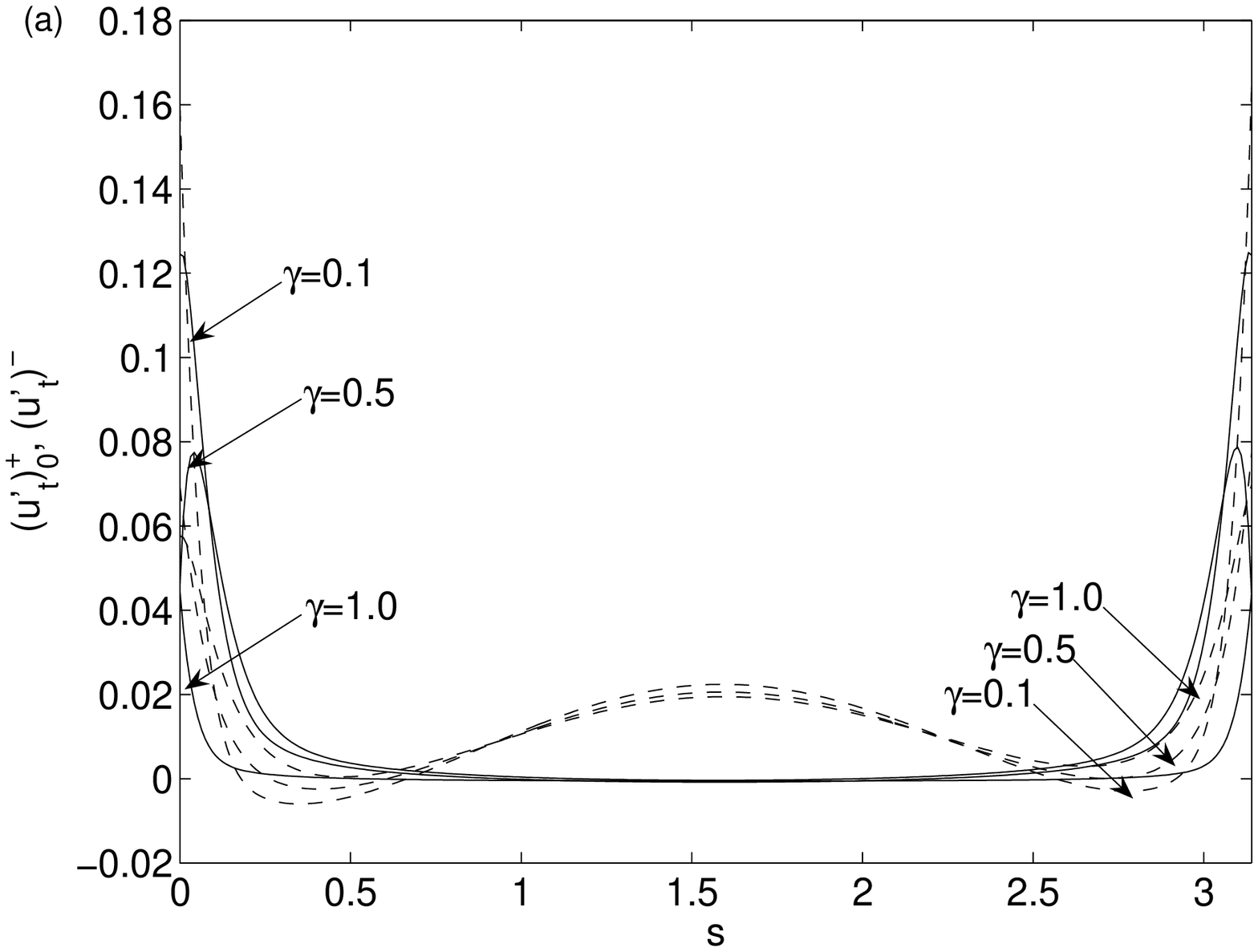} \hspace{1.0cm} \includegraphics{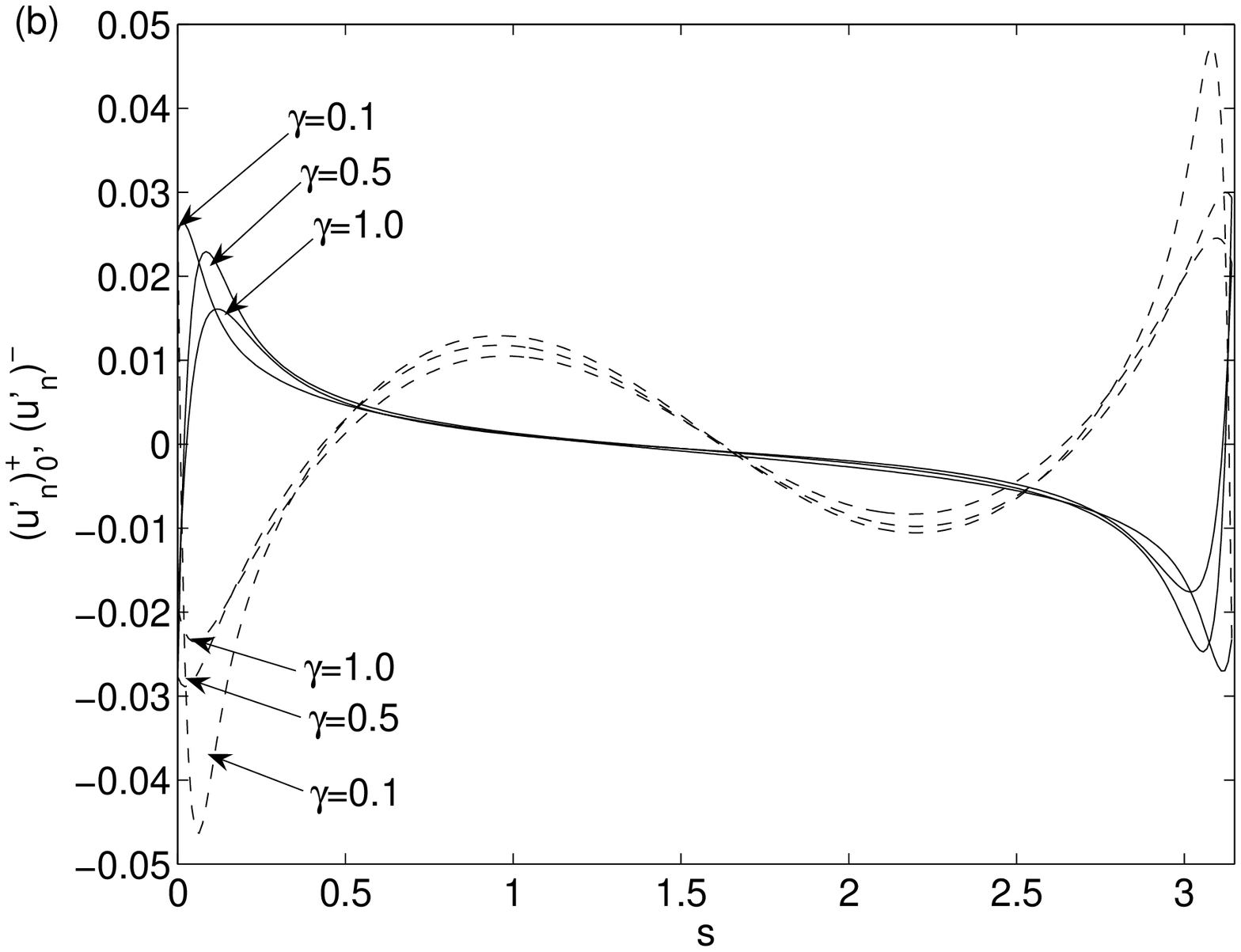}}
		\scalebox{0.3}{\includegraphics{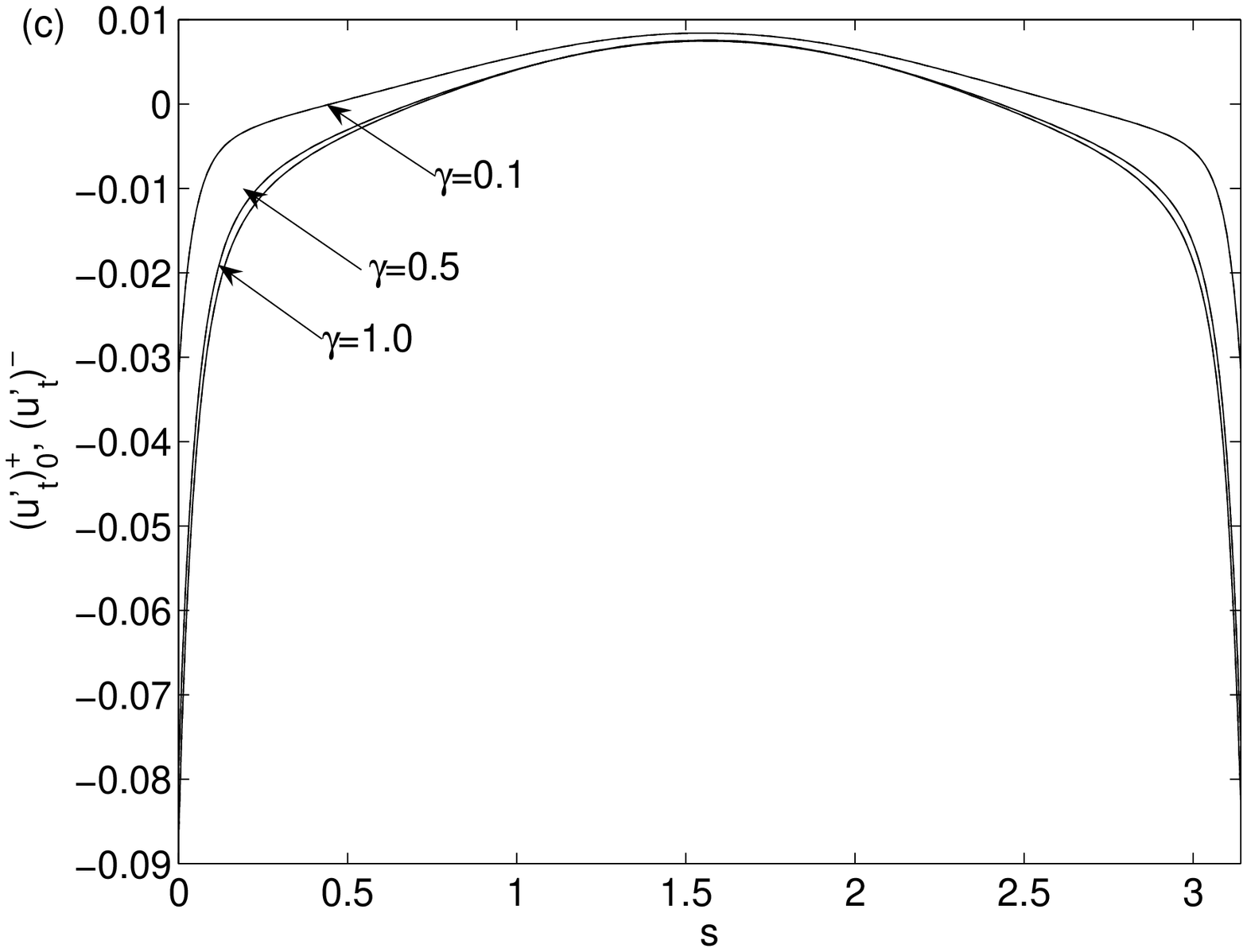} \hspace{1.0cm} \includegraphics{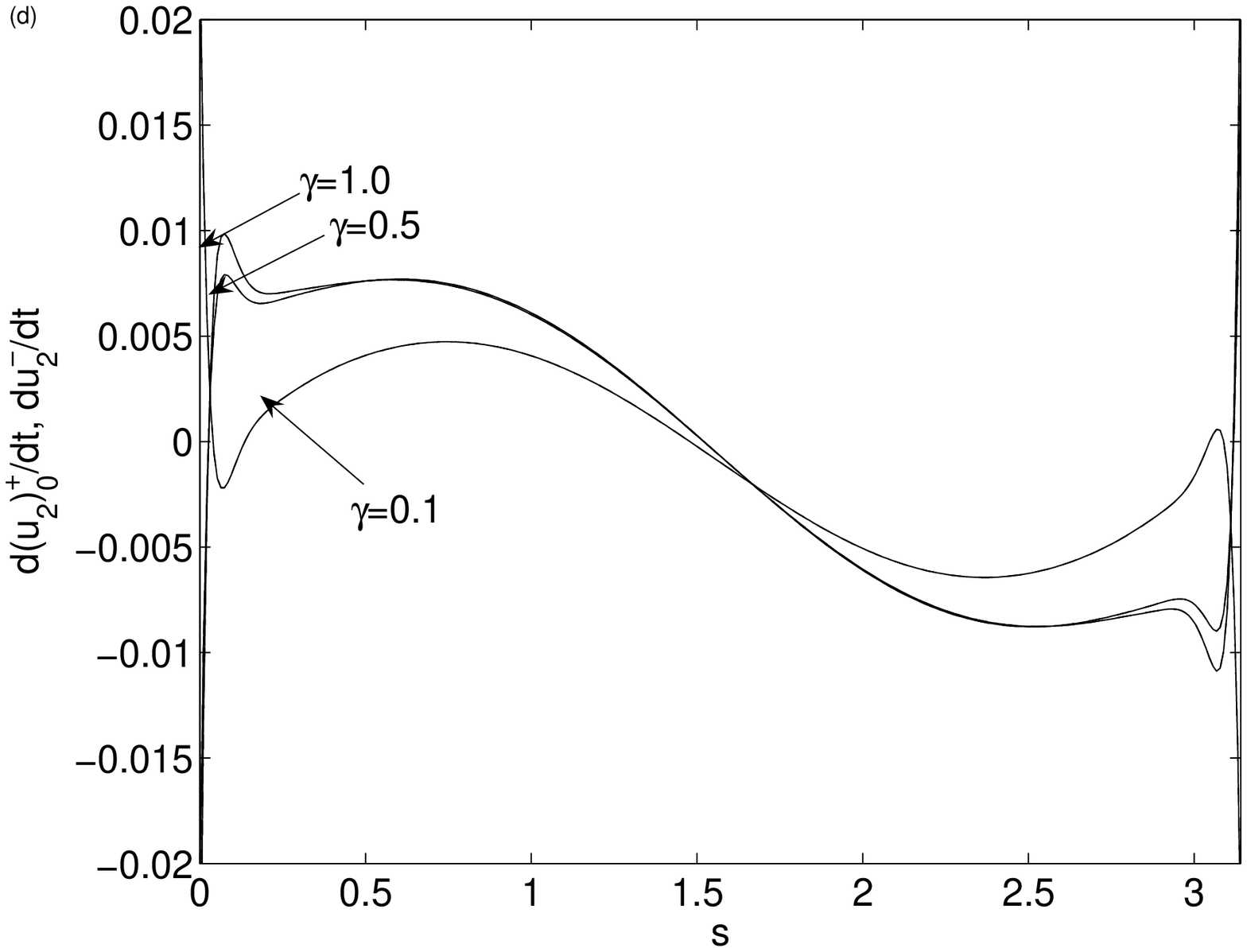}}
			\caption{Graphs of the derivatives of the displacements $(u'_t)_0^+$, $(u'_t)^-$, $(u'_n)_0^+$, $(u'_t)^-$ on the curves $L_0$ and $L$}
	\label{fig3}
\end{figure}

It can be seen on the examples that the numerical method described above produces the results with good accuracy. The fig. \ref{fig1} shows the graphs of the unknown functions $\R g'_0(s)$, $\I g'_0(s)$ computed on the interval $[0,l_0]$ for three values of the parameter $N$: $N=16$, $N=24$ and $N=30$. The solid (dashed) lines represent the graphs of the functions $\R g'_0(s)$ (or $\I g'_0(s)$). It can be seen that a good approximation of the results is achieved already for a relatively small value of $N=16$. The graphs above are for an upper semicircular fracture described by the parametric equation $t(s)=Re^{is}$, $s\in [0,\pi]$, on the boundary of a circular inclusion of a radius $R=1$. The mechanical parameters of the problem are as follows: $\mu=40$ GPa, $\nu=0.25$, $\mu_0=60$ GPa, $\nu_0=0.35$, $\gamma^+=\gamma^-=\gamma^i=0.1$, $\sigma^{\infty}_1=1$ MPa, $\sigma^{\infty}_2=0$ MPa, $\alpha=0$, $f_{1,2}(s)=0$.

\begin{figure}[ht]
	\centering
		\scalebox{0.3}{\includegraphics{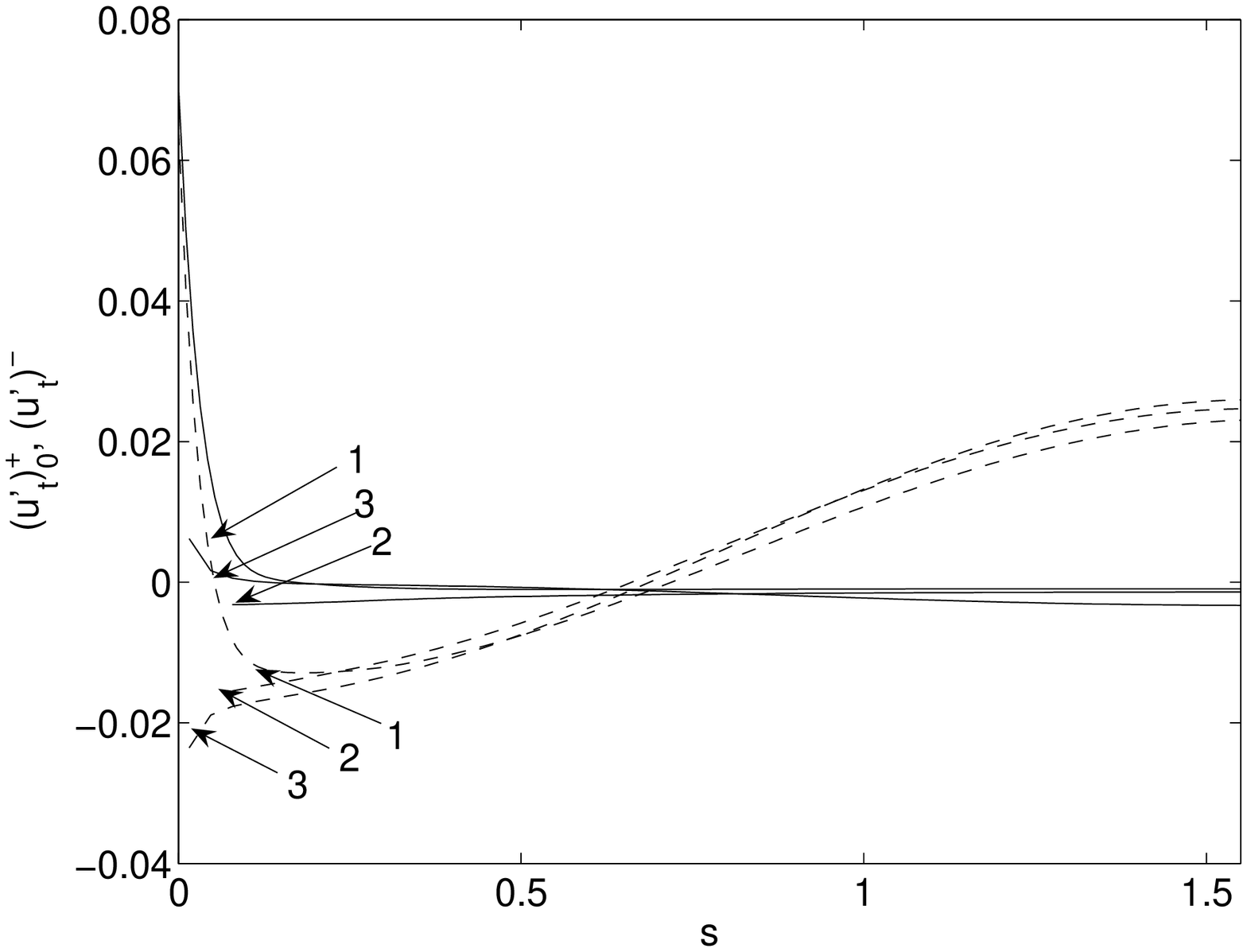} \hspace{1.0cm} \includegraphics{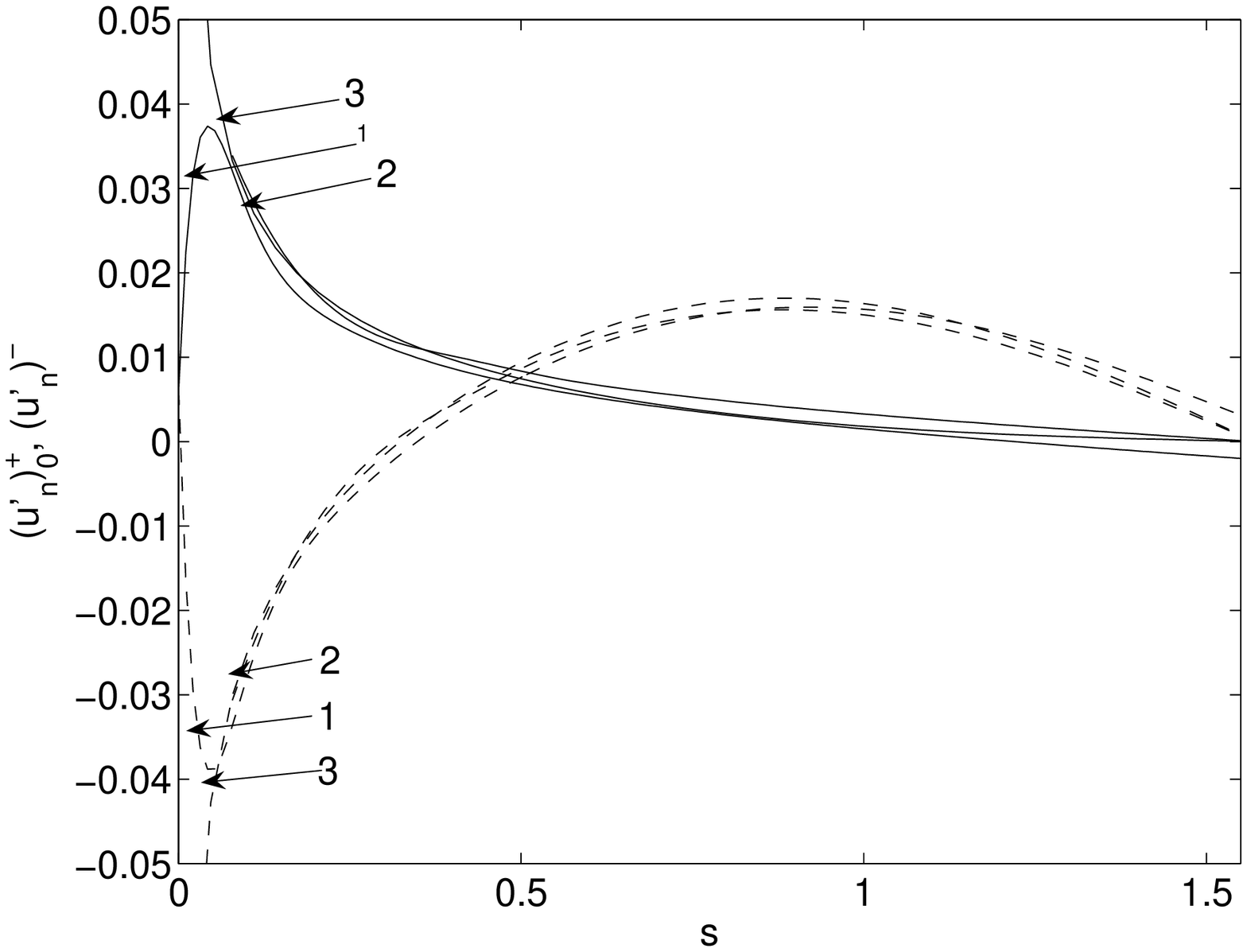}}
			\caption{Graphs of the derivatives of the displacements $(u'_t)_0^+$, $(u'_t)^-$, $(u'_n)_0^+$, $(u'_t)^-$ on the fracture surface $L_0$}
	\label{fig4}
\end{figure}

The fig. \ref{fig2} shows the graphs of the normal and shear stresses $\sigma^{+}_{n0}$, $\sigma_n^-$, $\tau^{+}_{n0}$, $\tau_n^-$ on the boundary of the semicircular interface fracture $L_0$ described by the parametric equation $t(s)=e^{is}$, $s\in [0, \pi]$ (fig. \ref{fig2}a, \ref{fig2}b) and on the connected part $L$ of the boundary of the circular inclusion  described by the parametric equation $t(s)=e^{i(s+\pi)}$, $s\in [0,\pi]$ (fig. \ref{fig2}c, \ref{fig2}d). The other parameters of the problem are as follows: $\mu=40$ GPa, $\nu=0.25$, $\mu_0=60$ GPa, $\nu_0=0.35$, $\gamma^+=\gamma^-=\gamma$, $\gamma^i=0$, $\sigma^{\infty}_1=1$ MPa, $\sigma^{\infty}_2=0$ MPa, $\alpha=0$, $f_{1,2}(s)=0$. The solid line on the graphs depicts the stresses $\sigma_{n0}^+$, $\tau_{n0}^+$ on the boundary of the inclusion $S_0$, the dashed line represents the stresses $\sigma_n^-$, $\tau_n^-$ on the boundary of the matrix $S$. The graphs are plotted for the following three values of the surface tension parameter $\gamma$: $\gamma=0.1$, $\gamma=0.5$ and $\gamma=1.0$. The graphs of the tangential and normal components $(u'_t)_0^+$, $(u'_t)^-$, $(u'_n)_0^+$, $(u'_t)^-$ of the derivatives of the displacements for the same values of the mechanical and elastic parameters are shown on the fig. \ref{fig3}.

\begin{figure}[ht]
	\centering
		\scalebox{0.3}{\includegraphics{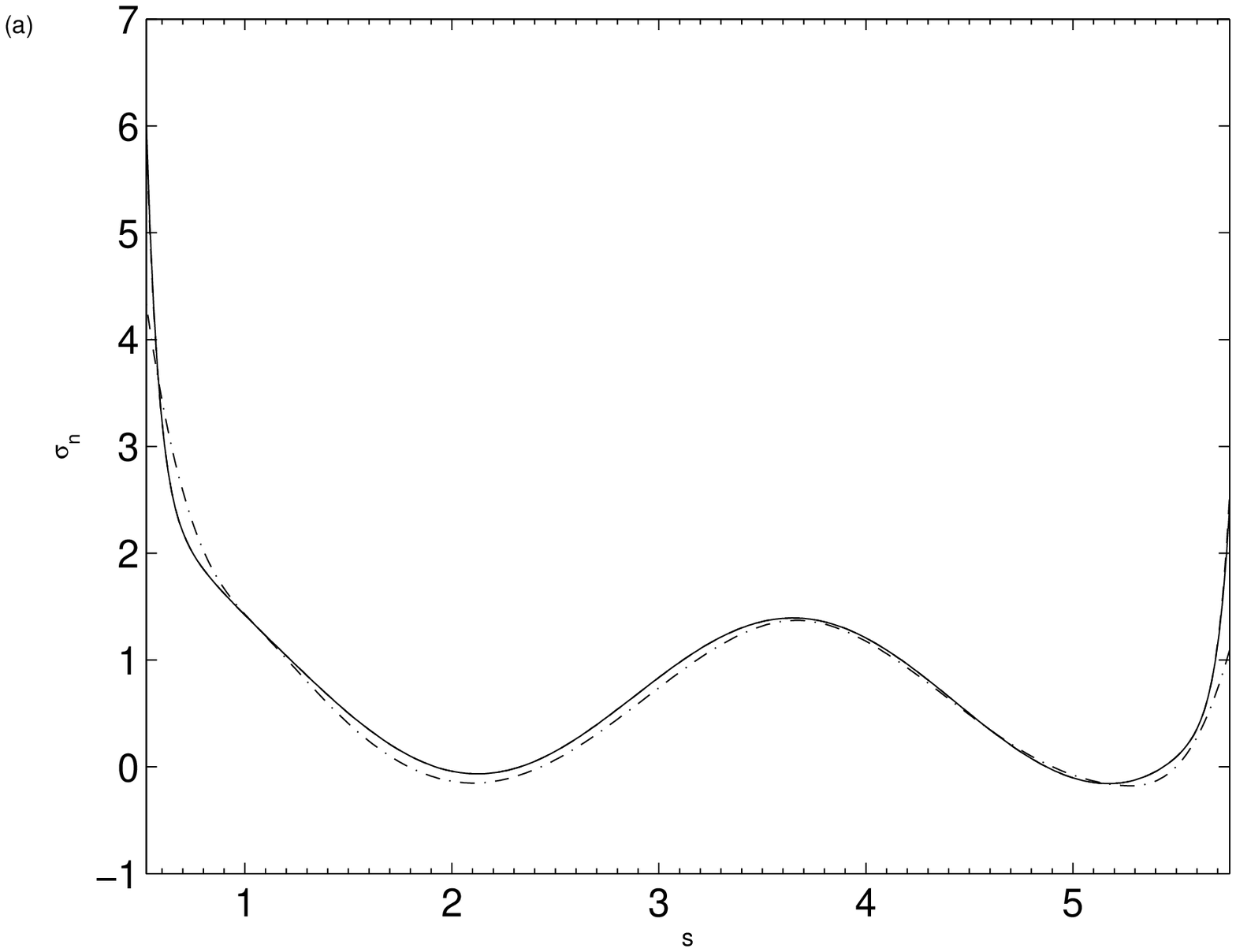} \hspace{1.0cm} \includegraphics{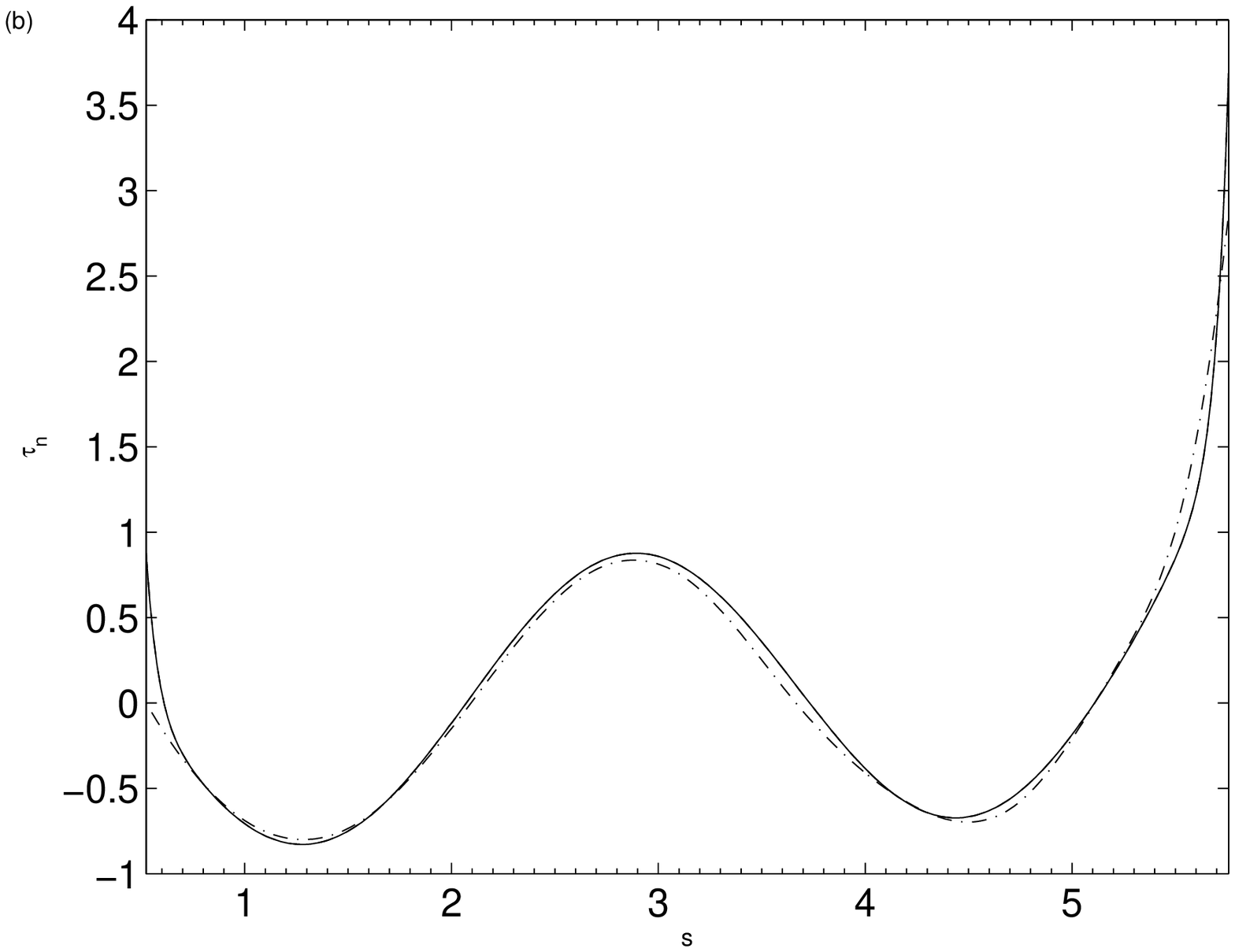}}
			\caption{Comparison of the stresses $\sigma_n$ and $\tau_n$ on the connected boundary of the circular inclusion with classical results}
	\label{fig5a}
\end{figure}

\begin{figure}[ht]
	\centering
		\scalebox{0.3}{\includegraphics{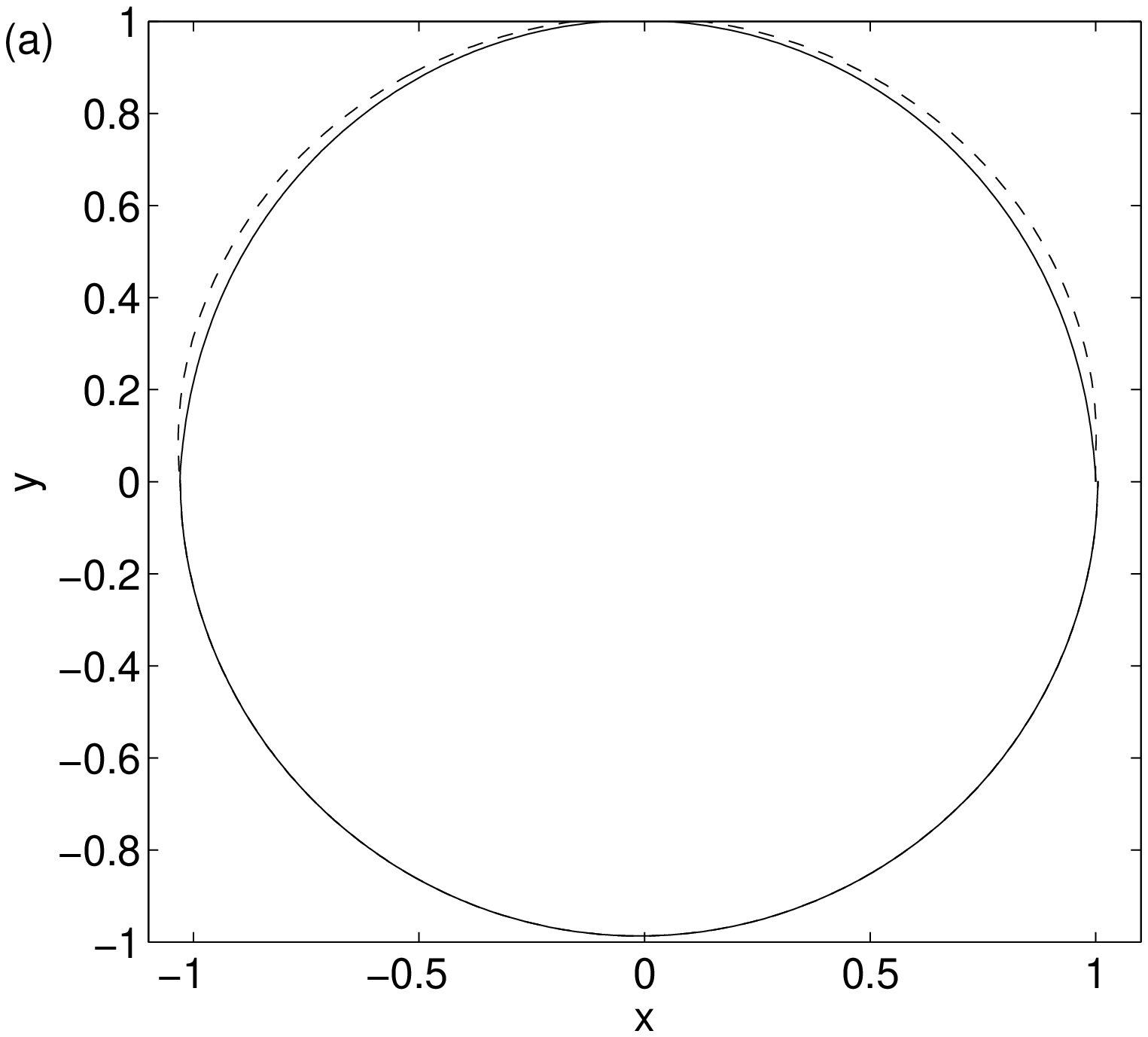} \hspace{1.0cm} \includegraphics{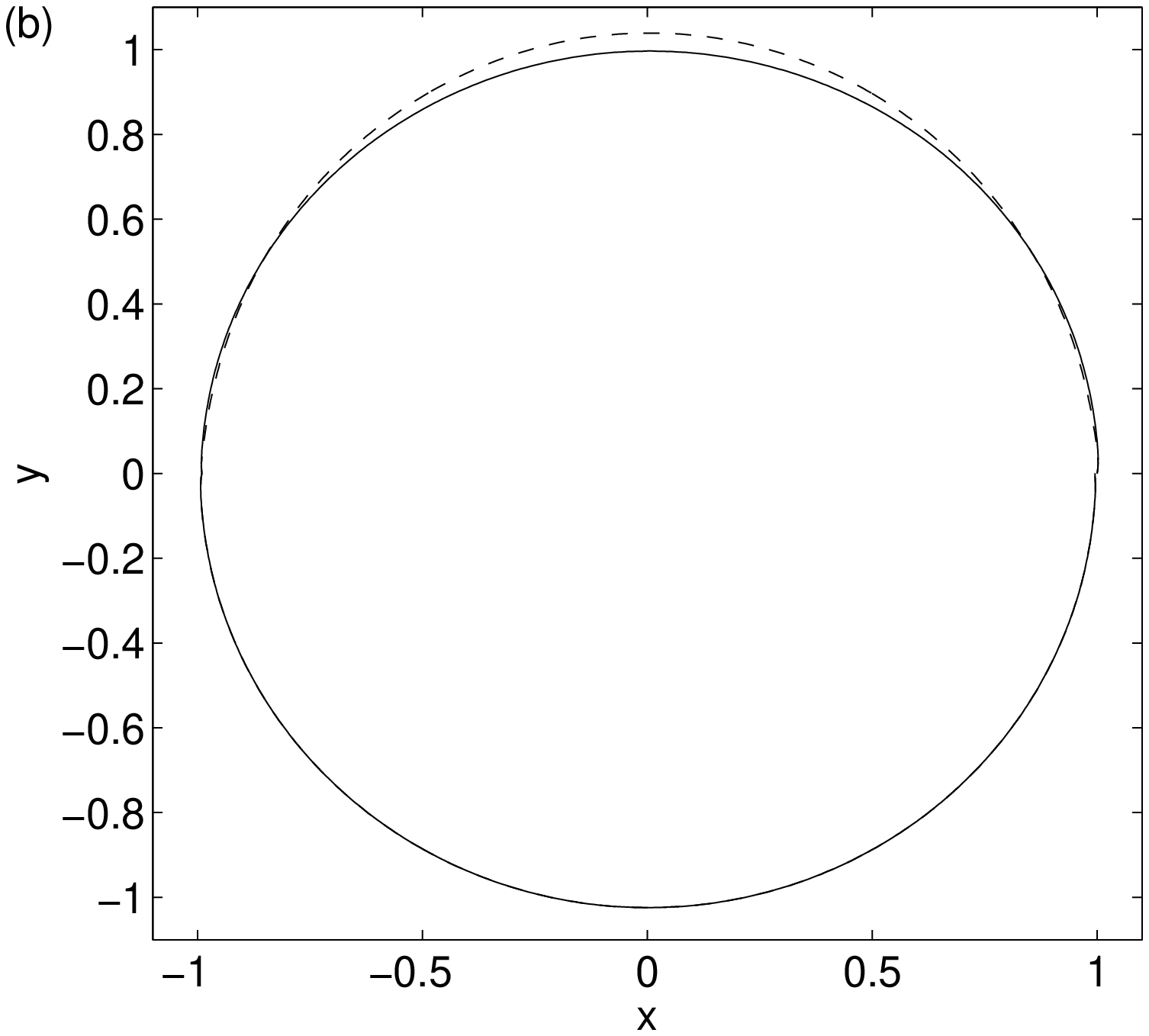}}
			\caption{Displacements of the boundary of the inclusion $S_0$ for (a) horizontal and (b) vertical remote stretching}
	\label{fig5}
\end{figure}

Comparison of the results obtained in this paper with previously known results is presented on the fig. \ref{fig4}. The solid lines on the fig. \ref{fig4} correspond to the graphs of the derivatives of the displacements $(u'_t)_0^+$, $(u'_n)_0^+$ on the right half of the semicircular fracture $L_0$ defined by the equation $t(s)=e^{is}$, $s\in [0,\pi/2]$, from the side of the inclusion $S_0$, the dashed lines correspond to the derivatives of the displacements $(u'_t)^-$, $(u'_t)^-$ on the same line from the side of the matrix $S$. The graphs on the left half of the line $L_0$ corresponding to the parameter values $s\in [\pi/2,\pi]$ can be obtained by symmetry. The lines denoted by the number ``1" correspond to the results for an interface fracture obtained using the methods presented in this paper, the lines denoted by ``2" correspond to the results for a non-interface fracture with a curvature-dependent surface tension obtained by similar methods in the paper \cite{Zemlyanova2012}, and finally, the lines denoted by ``3" correspond to the classical circular-arc-shaped interface fracture problem without surface tension on the fracture boundary which has been considered in different formulations in \cite{PrasadSimha2003}, \cite{SilZem2001}, \cite{Toya1974} among other works. The mechanical parameters for the problem are $\mu_0=\mu=40$ GPa, $\nu_0=\nu=0.25$, $\gamma^+=\gamma^-=0.01$, $\gamma^i=0.05$, $\sigma^{\infty}_1=1$ MPa, $\sigma^{\infty}_2=0$ MPa, $\alpha=0$, $f_{1,2}(s)=0$. In the case of the line ``2" the parameter $\gamma^i=0$, in the case of the line ``3" the parameters $\gamma^+=\gamma^-=\gamma^i=0$. Observe that because of the small differences in the initial values of the parameters it is possible to have small differences in the values of the derivatives of the displacements $(u'_t)_0^+$, $(u'_t)^-$, $(u'_n)_0^+$, $(u'_t)^-$. It can be seen from the fig. \ref{fig4} that the results are similar for all three considered problems except near the crack tip where the surface tension plays the major role and strongly affects the behavior of the stresses and the derivatives of displacements.

Comparison of the results obtained in this paper with the results for one interface fracture on the boundary of a circular inclusion presented in \cite{CrouchMogilevskaya2006} is given on the fig. \ref{fig5a}. The results are computed for a circular glass inclusion with elastic parameters $\mu_0=44.2$ GPa, $\nu_0=0.22$ inserted into an epoxy matrix with elastic parameters $\mu=2.39$ GPa, $\nu=0.35$. The crack occupies the region of the inclusion boundary defined by the polar angle $\theta\in [-\pi/6,\pi/6]$, while the connected part of the boundary is defined by the polar angle $\theta\in [\pi/6,11\pi/6]$. The stresses at infinity are given by $\sigma^{\infty}_1=1$ MPa, $\sigma^{\infty}_2=0$ MPa, $\alpha=\pi/6$. In the computations based on the results presented in the current paper it is assumed that the surface tension parameters $\gamma^+=\gamma^-=0.0001$, $\gamma^i=0$. The computations presented in \cite{CrouchMogilevskaya2006} are done for a classical problem with $\gamma^+=\gamma^-=\gamma^i=0$, hence, small differences in obtained results are possible. The fig. \ref{fig5a}a presents the graphs of the normal stresses $\sigma_n$ on the connected part of the boundary $L$, while the fig. \ref{fig5a}b presents the graphs of the shear stresses $\tau_n$ on $L$. The solid lines correspond to the results obtained in this paper, and the dashed lines correspond to the results obtained in \cite{CrouchMogilevskaya2006}. Good correspondence of the results can be observed.

The fig. \ref{fig5} shows the shape of the inclusion boundary and the interface fracture for two different types of the remote loading: (a) horizontal, $\sigma^{\infty}_1=1$ MPa, $\sigma^{\infty}_2=0$ MPa, $\alpha=0$, and (b) vertical $\sigma^{\infty}_1=1$ MPa, $\sigma^{\infty}_2=0$ MPa, $\alpha=\pi/2$. The displacements of the boundary $L_0\cup L$ have been scaled with a factor of two for better visibility. The mechanical parameters of the problem are as follows: $\mu=40$ GPa, $\nu=0.25$, $\mu_0=60$ GPa, $\nu_0=0.35$, $\gamma^+=\gamma^-=0.1$, $\gamma^i=0.05$, $f_{1,2}(s)=0$. The solid line represents the displacements of the boundary of the inclusion $S_0$ and the dashed line represents the displacements of the boundary of the matrix $S$. Observe, that since the contact conditions on the boundary of the fracture have not been imposed, it is possible for the boundaries of the crack to overlap.

\begin{figure}[ht]
	\centering
		\scalebox{0.3}{\includegraphics{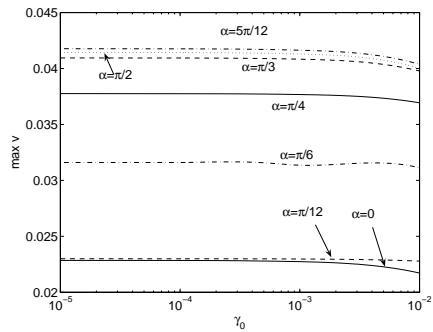}}
			\caption{Graphs of the maximal crack opening with respect to the surface tension parameter $\gamma_0$}
	\label{fig6}
\end{figure}

Finally, the fig. \ref{fig6} presents the graphs for the maximal crack opening defined as $\max v=\max_{t\in L_0}|d(u_1+iu_2)^+_0(t)/dt-d(u_1+iu_2)^-(t)/dt|$ for a circular crack with a parametric equation $t(s)=e^{is}$, $s\in [0,\pi]$. The material parameters of the inclusion and the matrix are given by $\mu=40$ GPa, $\nu=0.25$, $\mu_0=60$ GPa, $\nu_0=0.35$, $\gamma^+=\gamma^-=\gamma_0$, $\gamma^i=0$, $\sigma^{\infty}_1=1$ MPa, $\sigma^{\infty}_2=0$ MPa. The computations are made for different values of the parameter $\gamma_0$ and the angle $\alpha$ which specifies the direction of the principal stress $\sigma_1$ at infinity. It can be seen that the maximal crack opening decreases slightly as the value of $\gamma_0$ increases for all considered angles $\alpha$.

\section{Conclusions}
The paper is concerned with a study of an arbitrarily shaped interface fracture with a curvature-dependent surface tension acting on the boundary of the fracture. It has been shown that the fracture models which include surface effects on the fracture boundary can be efficient in eliminating integrable power crack tip singularities present in linear elasticity. Weaker logarithmic type singularities may still be present at the crack tips. This conclusion generalizes the results obtained in \cite{SendovaWalton2010}, \cite{Zemlyanova2013}, \cite{Zemlyanova2012}. Development of fracture models with surface effects which eliminate crack tip singularities completely remains an interesting open problem.  

\vspace{.1in}

\begin{acknowledgements}
The author gratefully acknowledges the support of the Simons Foundation through the AMS Simons Grant (2012-2014) and the Simons Collaboration Grant for Mathematicians (2014-2019).
\end{acknowledgements}

\end{document}